\documentclass[a4paper,12pt]{article}
\usepackage{a4wide}
\usepackage[ruled,vlined]{algorithm2e}
\usepackage{amsmath,amsfonts,amssymb,amsthm,url,textcomp,mathrsfs,tikz-cd}
\usepackage{anyfontsize}
\usepackage[font=scriptsize,labelfont=bf, justification=centering]{caption} 
\usepackage{csquotes}
\usepackage{enumitem}
\usepackage{floatrow}
\usepackage[a4paper, left=0.75in, right=0.75in, top=1.25in, bottom=1.75in]{geometry}
\usepackage{graphicx}
\usepackage[colorlinks=true, linkcolor=cyan]{hyperref}
\usepackage[capitalise]{cleveref}
\usepackage{indentfirst}
\usepackage{leftidx}
\usepackage{mathtools}
\usepackage{mdframed}
\usepackage[numbers]{natbib}
\usepackage{titling}
\usepackage{tikz}
\usetikzlibrary[positioning]
\usepackage{xcolor}

\hypersetup{citecolor=cyan}

\DeclareUrlCommand\url{\color{cyan}} 

\allowdisplaybreaks[1]

\makeatletter
\let\@fnsymbol\@arabic
\makeatother

\newtheorem{Cor}{Corollary}\numberwithin{Cor}{section}
\numberwithin{cor}{subsection}

\newtheorem{dfn}[Cor]{Definition}

\newtheorem{exm}[Cor]{Example}

\newtheorem{lem}[Cor]{Lemma}

\newtheorem{prp}[Cor]{Proposition}

\newtheorem{thm}[Cor]{Theorem}

\newtheorem*{thm*}{Theorem}
\newtheorem{manualtheoreminner}{Theorem}
\newenvironment{manualtheorem}[1]{%
  \renewcommand\themanualtheoreminner{#1}%
  \manualtheoreminner
}{\endmanualtheoreminner}

\newcommand*\samethanks[1][\value{footnote}]{\footnotemark[#1]}

\newcommand{\FF}[2]{\mathbb{F}_{{#1}^{#2}}}
\newcommand{\FP}[1]{\mathbb{F}_{#1}}

\newcommand{\N}{\text{N}}

\newcommand{\NN}{\mathbb{N}}
\newcommand{\PG}{\text{PG}}

\newcommand{\sgn}{\text{sgn}}

\newcommand{\cof}{\text{cof}}
\newcommand{\Spn}{\text{Span}}
\newcommand{\Tr}{\text{Tr}}

\begin{document}
\title{Permutation polynomials, projective polynomials, and bijections between \texorpdfstring{$\mu_{\frac{q^n-1}{q-1}}$}{} and \texorpdfstring{$\PG(n-1,q)$}{}}

\date{}

\author{ Tong Lin\samethanks[1] \and Qiang Wang\thanks{School of Mathematics and Statistics, Carleton University, 1125 Colonel By Drive, Ottawa ON K1S 5B6, Canada. The authors were supported by the Natural Sciences and Engineering Research Council of Canada (RGPIN-2023-04673).
\newline \emph{E-mail addresses}:
\href{mailto:tonglin4@cmail.carleton.ca}{\color{cyan}{tonglin4@cmail.carleton.ca}}\ (T. Lin), 
\href{mailto:wang@math.carleton.ca}{\color{cyan}{wang@math.carleton.ca}}\ (Q. Wang).
\newline \emph{Mathematics Subject Classification}: 11T06 (11B37).
\newline \emph{Keywords}: Finite fields; Permutation polynomials; Projective geometries; Projective polynomials.}
}

\maketitle

\vspace{-3em}

\abstract{
Using arbitrary bases for the finite field $\FF{q}{n}$ over $\FP{q}$, we obtain the generalized M\"obius transformations (GMTs), which are a class of bijections between the projective geometry $\PG(n-1,q)$ and the set of roots of unity $\mu_{\frac{q^n-1}{q-1}}\subseteq \FF{q}{n}$, where $n\geq 2$ is any integer. We also introduce a class of projective polynomials, using the properties of which we determine the inverses of the GMTs. Moreover, we study the roots of those projective polynomials, which lead to a three-way correspondence between partitions of $\FF{q}{n}^\ast,\mu_{\frac{q^n-1}{q-1}}$ and $\PG(n-1,q)$. Through this correspondence and the GMTs, we construct permutation polynomials of index $\frac{q^n-1}{q-1}$ over $\FF{q}{n}$. 
}

\section{Introduction}\label{sec1}
Throughout this paper, let $\FF{q}{n}$ be the finite field with $q^n$ elements, where $q$ is a prime power and $n\geq 2$. A polynomial $f\in \FF{q}{n}[x]$ \textit{permutes} $\FF{q}{n}$ if it is a bijection over $\FF{q}{n}$. Permutation polynomials over finite fields have long been a fundamental subject of great interest since they find their applications in a diversity of areas such as cryptography \cite{KEA22a,LLHQ21a,MP14a,SH98a}, coding theory \cite{DH13a,LC07a}, combinatorial designs \cite{DY06a} and block ciphers \cite{BEA12a}. Over the years, a vast number of classes of permutation polynomials have been researched; see, for example, \cite{AW07a, AGW09a, AGW11a,  GOOG24a, XH15a,  LQW18a, MRS22a, PUW23a, PL01a, QL23a,  QW19a, QW24a, YD11a,YD14a, ZYP16a} and the references therein. Among them, the study of permutation polynomials of the form $x^rh\left(x^s\right)$ has been popularized by numerous researchers using different techniques, many of which involve the multiplicative case of the AGW-Criterion, which is the following lemma.


\begin{lem}[\cite{AW07a,PL01a,WL91a,QW07a, MZ09}]\label{L11}
Let $r\in \NN$. Let $\ell,s\in \NN$ be such that $q^n-1=\ell s$. The polynomial $f(x)=x^rh\left(x^s\right)$ permutes $\FF{q}{n}$ if and only if $\gcd(r,s)=1$ and $g(x)=x^rh(x)^s$ permutes $\mu_\ell$, the set of all $\ell$-th roots of unity.
\end{lem}

In fact, every polynomial $f\in \FF{q}{n}[x]$ such that $f(0)=0$ can be uniquely written as $x^rh\left(x^s\right)$, where $r\geq 0$ is the \textit{vanishing order}  and $\ell =(q^n-1)/s$  is the \textit{index} of $f$.  The notion of index was introduced in \cite{AGW09a}. 

Let $r$ be any positive integer coprime to $q-1$. In this paper, we focus on constructing permutation polynomials of vanishing order $r$ and index $\frac{q^n-1}{q-1}$ over $\FF{q}{n}$. That is, we find polynomials $h\in \FF{q}{n}[x]$ such that $f(x)=x^rh\left(x^{q-1}\right)$ permutes $\FF{q}{n}$. Since $f(0)=0$ for each such $f$, it suffices to construct $h$ for which the corresponding $f$ permutes $\FF{q}{n}^\ast$. For this purpose, we apply \cref{L11}. Through the mapping $x^{q-1}$, we project $\FF{q}{n}^\ast$ onto $\mu_{\frac{q^n-1}{q-1}}$, which is of smaller cardinality, and find $h$ such that $g(x)=x^rh(x)^{q-1}$ permutes $\mu_{\frac{q^n-1}{q-1}}$. To study the permutational behavior of $g$, we may follow the same ideas as above and project $\mu_{\frac{q^n-1}{q-1}}$ onto a set $S$ via properly chosen bijections from $\mu_{\frac{q^n-1}{q-1}}$ to $S$. For example, $S$ can be the projective geometry
$$\PG(n-1,q)=\bigl(\FP{q}^n\setminus \bigl\{\vec{0}\bigr\}\bigr)\slash \sim,$$
which is the quotient set of $\FP{q}^n\setminus \bigl\{\vec{0}\bigr\}$ by the equivalence relation $\sim$, under which $\vec{x}\sim \vec{y}$ when $\vec{x}=\lambda\vec{y}$ for some $\lambda\in \FP{q}^\ast$. Counting the number of equivalence classes yields that $\bigl|\PG(n-1,q)\bigr|=\frac{q^n-1}{q-1}=\Bigl|\mu_{\frac{q^n-1}{q-1}}\Bigr|$. We denote the equivalence class in $\PG(n-1,q)$ represented by a vector $\vec{x}=(x_0,\dots,x_{n-1})\in \FP{q}^n\setminus \bigl\{\vec{0}\bigr\}$ by $(x_0:\dots:x_{n-1})$.

\begin{figure}[hbt!]
\begin{tikzcd}[row sep=4ex, column sep=14ex]
{\FF{q}{n}^\ast} \arrow{r}{f(x)=x^rh\left(x^{q-1}\right)} \arrow{d}{x^{q-1}} & {\FF{q}{n}^\ast} \arrow{d}{x^{q-1}} \\
{\mu_{\frac{q^n-1}{q-1}}} \arrow{r}{g(x)=x^rh(x)^{q-1}} \arrow{d}{\eta_1^{-1}} & {\mu_{\frac{q^n-1}{q-1}}} \\
{\PG(n-1,q)} \arrow{r}{\overline{g}} & {\PG(n-1,q)} \arrow{u}[right]{\eta_2} 
\end{tikzcd}{\caption{The multiplicative AGW criterion with generic bijections $\eta_1$ and $\eta_2$}\label{Fig1}}
\end{figure}

As \cref{Fig1} shows, for each bijection $\overline{g}$ over $\PG(n-1,q)$ and each pair of bijections $\eta_1,\eta_2$ from $\PG(n-1,q)$ to $\mu_{\frac{q^n-1}{q-1}}$, the mapping $\eta_2\circ \overline{g}\circ \eta_1^{-1}$ is a bijection over $\mu_{\frac{q^n-1}{q-1}}$. If there exists an $h\in \FF{q}{n}[x]$ such that
\begin{equation}\label{Eq001}
x^rh(x)^{q-1}=\eta_2\circ \overline{g}\circ \eta_1^{-1},
\end{equation}
then the polynomial $f(x)=x^rh\left(x^{q-1}\right)$ permutes $\FF{q}{n}$ as per \cref{L11}. Therefore, our main goal in this paper becomes finding polynomials $h\in \FF{q}{n}[x]$ which satisfy \cref{Eq001} for bijections $\overline{g}$ over $\PG(n-1,q)$ and some bijections $\eta_1,\eta_2$ from $\PG(n-1,q)$ to $\mu_{\frac{q^n-1}{q-1}}$.

When $n=2$, a common choice for $\eta_1,\eta_2$ is the so-called M\"obius transformation. More precisely, for $i=1,2$, we pick $u_i\in \mu_{q+1}\setminus \FP{q}$, and let
$$\eta_i(x)=\frac{x+u_i^q}{x+u_i},$$
Since $\gcd(r,q-1)=1$, $f(x)=x^rh\left(x^{q-1}\right)$ permutes $\FF{q}{2}$ if and only if $\eta_2^{-1}\circ x^rh(x)^{q-1}\circ \eta_1$ permutes $\PG(1,q)\cong \FP{q}\sqcup \{\infty\}$, where $\infty$ denotes a `point at infinity' adjoined to the `line' $\FP{q}$. Examples of this technique can be found in \cite{GOOG24a,XH18a,LH17a,PUW23a,TZ19a,TZLH18a}.

Similarly, when $n=3$, \cref{Fig1} shows that one can construct $h\in \FF{q}{3}[x]$ for which $f(x)=x^rh\left(x^{q-1}\right)$ permutes $\FF{q}{3}$ by projecting $\mu_{q^2+q+1}$ onto $\PG(2,q)$. In this case, instead of using the usual M\"obius transformations, Qu and Li \cite{QL23a} introduced the following mappings. Let $W=\Bigl\{\omega,\omega^q,\omega^{q^2}\Bigr\}$ be a basis $\FF{q}{3}$ over $\FP{q}$ for some $\omega\in \FF{q}{3}^\ast$, and define $\psi_W:\PG(2,q)\to \mu_{q^2+q+1}$, where
\begin{equation}\label{Eq002}
\psi_W(x_0:x_1:x_2)=\frac{x_0\omega^q+x_1\omega^{q^2}+x_2\omega}{x_0\omega+x_1\omega^q+x_2\omega^{q^2}}.
\end{equation}
Let $\alpha=\omega^{q+1}-\omega^{2q^2}$ and define, for $i=0,1,2$,
\begin{equation}\label{Eq003}
T_{W,i}(x)=\alpha^{q^{i+1}}+\alpha^{q^{i+2}}x+\alpha^{q^i}x^{q+1},
\end{equation}
Then for each $x\in \mu_{q^2+q+1}$,
\begin{equation}\label{Eq004}
\psi_W^{-1}(x)=
\begin{cases}
(1:0:0) & \text{if }x=\omega^{q-1};\\
\biggl(\frac{T_{W,0}(x)}{T_{W,1}(x)}:1:0\biggr) & \text{if }T_{W,1}(x)\neq 0\text{ and }T_{W,2}(x)=0;\\
\biggl(\frac{T_{W,0}(x)}{T_{W,2}(x)}:\frac{T_{W,1}(x)}{T_{W,2}(x)}:1\biggr) & \text{if }T_{W,2}(x)\neq 0.
\end{cases}
\end{equation}
In this paper, we generalize Qu and Li's ideas in order to deal with the case where $n\geq 2$ is an arbitrary integer. We summarize our first main result into the following theorem.

\begin{manualtheorem}{I}\label{MR1}
Let $W=\{\omega_j:0\leq j\leq n-1\}$ be any basis for $\FF{q}{n}$ over $\FP{q}$. Define $\psi_W:\PG(n-1,q)\to \mu_{\frac{q^n-1}{q-1}}$, where
$$\psi_W(x_0:\dots:x_{n-1})
=\Biggl(\sum_{j=0}^{n-1}x_j\omega_j\Biggr)^{q-1}
=\frac{\displaystyle{\sum_{j=0}^{n-1}x_j\omega_j^q}}{\displaystyle{\sum_{j=0}^{n-1}x_j\omega_j}}.$$
Then $\psi_W$ is a bijection.
Moreover, for $0\leq i\leq n-1$, define
$$T_{W,i}(x)=\sum_{k=0}^{n-1}(-1)^{k(n-1)}c_i^{q^k}x^{\frac{q^k-1}{q-1}},$$
in which the constant term $c_i=T_{W,i}(0)$ is the $(i,0)$-cofactor of the matrix
$$M_W=\begin{bmatrix}
\omega_0     & \omega_0^q     & \dots  & \omega_0^{q^{n-1}} \\
\omega_1     & \omega_1^q     & \dots  & \omega_1^{q^{n-1}} \\
\vdots       & \vdots         & \ddots & \vdots \\
\omega_{n-1} & \omega_{n-1}^q & \dots  & \omega_{n-1}^{q^{n-1}} \\
\end{bmatrix}.$$
Then
\[
\resizebox{0.8\hsize}{!}{%
$\psi_W^{-1}(x)
=\begin{cases}
(1:0:\dots:0) & \text{ if }x=\omega_0^{q-1}; \\
\displaystyle{\biggl(\frac{T_{W,0}(x)}{T_{W,1}(x)}:1:0:\dots:0\biggr)} & \text{ if }x=(x_0\omega_0+\omega_1)^{q-1}\text{ for some }x_0\in \FP{q};\\
\vdots &\ \vdots \\
\displaystyle{\biggl(\frac{T_{W,0}(x)}{T_{W,n-1}(x)}:\dots:\frac{T_{W,n-2}(x)}{T_{W,n-1}(x)}:1\biggr)} & \begin{array}{l}\text{if }x=(x_0\omega_0+\dots+x_{n-2}\omega_{n-2}+\omega_{n-1})^{q-1} \\ \text{for some }x_0,\dots,x_{n-2}\in \FP{q}.\end{array}
\end{cases}$%
}
\]
\end{manualtheorem}

\begin{figure}[hbt!]
\begin{tikzcd}[row sep=4ex, column sep=14ex]
{\FF{q}{n}^\ast} \arrow{r}{f(x)=x^rh\left(x^{q-1}\right)} \arrow{d}{x^{q-1}} & {\FF{q}{n}^\ast} \arrow{d}{x^{q-1}} \\
{\mu_{\frac{q^n-1}{q-1}}} \arrow{r}{g(x)=x^rh(x)^{q-1}} \arrow{d}{\psi_W^{-1}} & {\mu_{\frac{q^n-1}{q-1}}} \\
{\PG(n-1,q)} \arrow{r}{\overline{g}} & {\PG(n-1,q)} \arrow{u}[right]{\psi_Y} 
\end{tikzcd}{\caption{The multiplicative AGW criterion with generalized M\"{o}bius transformations}\label{Fig2}}
\end{figure}

As \cref{Fig2} shows, if $W,Y$ are any bases for $\FF{q}{n}$ over $\FP{q}$, the GMT $\psi_W$, the inverse of the GMT $\psi_Y$ together with any bijection $\overline{g}$ over $\PG(n-1,q)$ allow us to construct bijections over $\mu_{\frac{q^n-1}{q-1}}$, which then induce permutations over $\FF{q}{n}$ via \cref{L11}.

When $n=2$, and $W=\{1,u\}$, where $u\in \mu_{q+1}\setminus \{\pm 1\}$, it is easy to check that a usual M\"obius transformation and its inverse satisfy the formulas given in \cref{MR1}. When $n=3$ and $W=\left\{\omega,\omega^q,\omega^{q^2}\right\}$ for some $\omega\in \FF{q}{3}^\ast$, \cref{MR1} recovers Qu and Li's results in \cite{QL23a}.

The polynomials $T_{W,i}\ (0\leq i\leq n-1)$ are examples of \textit{projective polynomials}, which are of the form $\displaystyle{\sum_{k=0}^{n-1}c_kx^{\frac{q^k-1}{q-1}}}$. They were introduced by Abhyankar \cite{SA97a}. As a side note, our scope of understanding of projective polynomials is rather limited. Many known results about projective polynomials concern those with very few terms. For example, projective trinomials of the form $x^{p^\ell+1}+x+a$ over $\mathbb{F}_{p^n}$ were studied by Berlekemp et al. when $p=2$ and $\ell=1$ \cite{BRS66a}. Their works were subsequently generalized by Helleseth and Kolosha when $p=2$ and $\ell$ is arbitrary \cite{HK10a}, and Bluher when $p$ and $\ell$ are both arbitrary \cite{AB04a}. In full generality, Kim et al. characterized the roots the above-mentioned trinomials \cite{KCM21a}. We also refer the readers to \cite{MS19a}, which contains more references on projective polynomials and their applications.\\
\\
Using the roots of $T_{W,i}\ (0\leq i\leq n-1)$, our next main result demonstrates a three-way correspondence between partitions of $\FF{q}{n}^\ast$, $\mu_{\frac{q^n-1}{q-1}}$ and $\PG(n-1,q)$. This correspondence plays a key role in our construction of permutation polynomials of index $\frac{q^n-1}{q-1}$ over $\FF{q}{n}$.

\begin{manualtheorem}{II}\label{MR2}
Let $W=\{\omega_j:0\leq j\leq n-1\}$ be any basis for $\FF{q}{n}$ over $\FP{q}$, and let $B=\{\beta_j:0\leq j\leq n-1\}$ be the dual basis of $W$. That is, $\Tr_q^{q^n}(\beta_i\omega_j)$ is $1$ if $i=j$ and $0$ otherwise. Let $T_{W,i}\ (0\leq i\leq n-1)$ be as defined in \cref{MR1}. For $0\leq j\leq n-1$, define
\begin{align*}
S_{W,j} & =\Bigl\{y\in \FF{q}{n}^\ast:\Tr_q^{q^n}(\beta_jy)\neq 0\text{ and }\Tr_q^{q^n}(\beta_iy) =0 \text{ for all }i\geq j+1\Bigr\}; \\
Z_{W,j} & =\bigl\{x\in \mu_{\frac{q^n-1}{q-1}}:T_{W,j}(x)\neq 0\text{ and }T_{W,i}(x)=0\text{ for all }i\geq j+1\bigr\}; \\
\overline{C_{q,n,j}} & =\biggl\{\biggl(\frac{Tr_q^{q^n}(\beta_0x)}{Tr_q^{q^n}(\beta_jx)}:\dots:\frac{Tr_q^{q^n}(\beta_{j-1}x)}{Tr_q^{q^n}(\beta_jx)}:1:0:\dots:0\biggr):x\in S_{W,j}\biggr\}.
\end{align*}
Then we have the following partitions.
\begin{align*}
\FF{q}{n}^\ast=\bigsqcup_{j=0}^{n-1}S_{W,j};\quad
\mu_{\frac{q^n-1}{q-1}}=\bigsqcup_{j=0}^{n-1}Z_{W,j};\quad
\PG(n-1,q)=\bigsqcup_{j=0}^{n-1}\overline{C_{q,n,j}}.
\end{align*}
Moreover, for $0\leq j\leq n-1$,
$$Z_{W,j}=\psi_W\bigl(\overline{C_{q,n,j}}\bigr)=S_{W,j}^{q-1}.$$
\end{manualtheorem}

The three partitions in \cref{MR2} allow us to construct permutational polynomials of the form $f(x)=x^rh\left(x^{q-1}\right)$ over $\FF{q}{n}^\ast$ by building $h$ piecewise on $S_{W,j}\ (0\leq j\leq n-1)$. More precisely, let $W,Y$ be two bases for $\FF{q}{n}$ over $\FP{q}$, and let $B=\left\{\beta_j:0\leq j\leq n-1\right\}$ be the dual basis of $W$. Define $\psi_W,\psi_Y$ by \cref{MR1}, and let $\overline{g}=(g_0:\dots:g_{n-1})$ be any bijection over $\PG(n-1,q)$, where $g_0,\dots,g_{n-1}\in \FP{q}[x_0,\dots,x_{n-1}]$. Assume that
$$h(x) =\begin{cases}
\widetilde{h}(x) & \text{if }x\in \FF{q}{n}\setminus \mu_{\frac{q^n-1}{q-1}}; \\
h_j(x) & \text{if }x\in Z_{W,j}\ (0\leq j\leq n-1),
\end{cases}$$
where $\widetilde{h}$ is an arbitrary polynomial over $\FF{q}{n}$, and on each $Z_{W,j}$,
$$x^rh_j(x)^{q-1}=\psi_Y\circ \overline{g}\circ \psi_W^{-1}(x).$$
Also, define $g:\mu_{\frac{q^n-1}{q-1}}\to \mu_{\frac{q^n-1}{q-1}}$, where $g(x)=x^rh_j(x)^{q-1}$ on each $Z_{W,j}$. For each $j$, let
$$f_j(x) =x^rh_j\bigl(x^{q-1}\bigr),$$
and define
$$f(x)=\begin{cases}
0   & \text{if }x=0; \\
f_j(x) & \text{if }x\in S_{W,j}\ (0\leq j\leq n-1).
\end{cases}$$
Clearly, $f(x)=x^rh\bigl(x^{q-1}\bigr)$ and $g=\psi_Y\circ \overline{g}\circ \psi_W^{-1}$. In particular, $g$ permutes $\mu_{\frac{q^n-1}{q-1}}$. So by \cref{L11}, $f$ is permutes $\FF{q}{n}$. Based on these ideas, we present our third main result, which enables us to construct \textit{homogeneous} permutational polynomials of index $\frac{q^n-1}{q-1}$ over $\FF{q}{n}$. We recall that a mapping $f\in \FF{q}{n}[x]$ is called homogeneous of degree $r$, or simply $r$-homogeneous for some $r\in \NN\sqcup \{0\}$ if $f(\lambda x)=\lambda^rf(x)$ for all $\lambda\in \FP{q}$ and $x\in \FF{q}{n}$.

\begin{manualtheorem}{III}\label{MR3}
Let $W$ and $Y=\left\{\gamma_j:0\leq j\leq n-1\right\}$ be any two bases for $\FF{q}{n}$ over $\FP{q}$. Let $B=\left\{\beta_j:0\leq j\leq n-1\right\}$ be the dual basis of $W$. Let $\overline{g}=(g_0:\dots:g_{n-1})$ be any bijection over $\PG(n-1,q)$, where $g_0,\dots,g_{n-1}:\FP{q}^n\to \FP{q}$ are multivariate polynomials. Let $r\in \NN$ be coprime to $q-1$. For $0\leq j\leq n-1$ and each $x\in \FF{q}{n}^\ast$ such that 
$$j=\max\bigl\{0\leq i\leq n-1:\Tr_q^{q^n}(\beta_ix)\neq 0\bigr\},$$
we define
$$f_j(x)=\Tr_q^{q^n}(\beta_jx)^r\sum_{k=0}^{n-1}g_k\Biggl(\frac{\Tr_q^{q^n}(\beta_0x)}{\Tr_q^{q^n}(\beta_jx)},\dots,\frac{\Tr_q^{q^n}(\beta_{j-1}x)}{\Tr_q^{q^n}(\beta_jx)},1,0,\dots,0\Biggr)\gamma_k.$$
Then  $f:\FF{q}{n}\to \FF{q}{n}$, where
$$f(x)=\begin{cases}
0 & \text{if }x=0; \\
f_j(x) & \text{if }j=\max\bigl\{0\leq i\leq n-1:\Tr_q^{q^n}(\beta_ix)\neq 0\bigr\}\quad (0\leq j\leq n-1),
\end{cases}$$
is a rational mapping permuting $\FF{q}{n}$. In particular, if $r\geq \max\{\deg(g_k):0\leq k\leq n-1\}$, then $f$ is a permutation polynomial of index $\frac{q^n-1}{q-1}$ over $\FF{q}{n}$.
\end{manualtheorem}

From the definition of $f$ in \cref{MR3}, it is easy to see that for all $\lambda\in \FP{q}$ and $x\in \FF{q}{n}$. $f_j(\lambda x)=\lambda^rf_j(x)\ (0\leq j\leq n-1)$.
Hence, $f$ is $r$-homogeneous.\\
\\
In order to apply \cref{MR3}, we need to find bijections $\overline{g}$ over $\PG(n-1,q)$. In fact, they can be built via vector-valued mappings over $\FP{q}^n$ from the following criterion. Let $\overline{E_i},E_i\ \bigl(0\leq i\leq \frac{q^n-1}{q-1}-1\bigr)$ denote all equivalence classes in $\PG(n-1,q)$ and their underlying sets of representatives in $\FP{q}^n\setminus \bigl\{\vec{0}\bigr\}$, and let $\widetilde{g}=(g_0,\dots,g_{n-1})$. Then $\overline{g}=(g_0:\dots:g_{n-1})$ is a bijection over $\PG(n-1,q)$ if and only if there exists a permutation $\sigma$ in the symmetric group $S_{\frac{q^n-1}{q-1}}$ such that $\widetilde{g}(E_i)=E_{\sigma(i)}$ for $0\leq i\leq \frac{q^n-1}{q-1}-1$. In particular, $\widetilde{g}$ need not be a bijection over $\FP{q}^n$.

We give a concrete example of \cref{MR3} here, and more can be found in \cref{sec3}. Let $H$ be a polynomial over $\FF{q}{2}$ which permutes $\FP{q}$. Let $g_0,g_1:\FP{q}^2\to \FP{q}$ be such that $(g_0,g_1)\bigl(\vec{0}\bigr)=\vec{0}$, and for each $(x_0,x_1)\in \FP{q}^2\setminus \bigl\{\vec{0}\bigr\}$,
$$(g_0,g_1)(x_0,x_1)=
\begin{cases}
\Bigl(H\Bigl(\frac{x_0}{x_1}\Bigr),1\Bigr) & \text{if }x_1\neq 0; \\
(1,0) & \text{if }x_1=0.
\end{cases}$$
Then the corresponding mapping $\overline{g}=(g_0:g_1)$ over $\PG(n-1,q)$ is given by
$$\begin{cases}
\overline{g}(x_0:1)=(H(x_0):1) & \text{for all }x_0\in \FP{q}; \\
\overline{g}(1:0)=(1:0).
\end{cases}$$
Since $\overline{g}$ permutes $\PG(n-1,q)$, plugging $g_0,g_1$ in \cref{MR3} leads to the following class of permutation polynomials of index $q+1$ over $\FF{q}{2}$.

\begin{manualtheorem}{IV}\label{MR4}
Let $r,d\in \NN$ be such that $r\geq d$ and $\gcd(r,q-1)=1$. Let $H$ be a polynomial over $\FF{q}{2}$ of degree $d$ which permutes $\FP{q}$. Let $a,b\in \FF{q}{2}^\ast$ and $u,v\in \mu_{q+1}$ be such that $u\neq v$ and $\bigl(\frac{u}{v}\bigr)^d\neq \bigl(\frac{a}{b}\bigr)^{q-1}$. Let $w\in \FF{q}{2}^\ast$ be such that $w^{q-1}=\frac{u}{v}$. Then $F$ permutes $\FF{q}{2}$, where
$$F(x)=\begin{cases}
aw^{r-d}(v-u)^rx^r & \text{if }x^q-vx=0; \\
\bigl(x^q-vx\bigr)^r\left(aw^{-d}H\left(w\left(\frac{x^q-ux}{x^q-vx}\right)\right)+b\right) & \text{if }x^q-vx\neq 0.
\end{cases}$$
In particular, if $\gcd(d,q-1)=1$ and $H(x)=x^d$, then $\widetilde{F}$ permutes $\FF{q}{2}$, where
$$\widetilde{F}(x)=\begin{cases}
aw^{r-d}(v-u)^rx^r & \text{if }x^q-vx=0; \\
a\bigl(x^q-vx\bigr)^{r-d}\bigl(x^q-ux\bigr)^d+b\bigl(x^q-vx\bigr)^r & \text{if }x^q-vx\neq 0.
\end{cases}$$
\end{manualtheorem}

The rest of this paper is organized as follows. In \cref{sec2}, we generalize the M\"obius transformations and prove \cref{MR1}. In \cref{sec3}, we prove \cref{MR2} and \cref{MR3}. Besides, we verify \cref{MR4} and give more examples of homogeneous permutation polynomials of index $\frac{q^n-1}{q-1}$ over $\FF{q}{n}$. In \cref{sec4}, we compare the generalized M\"obius transformations to a known class of bijections from $\PG(n-1,q)$ to $\bigl\{\alpha^i:0\leq i\leq \frac{q^n-1}{q-1}-1\bigr\}$ for subprimitive roots $\alpha$ in $\FF{q}{n}$ \cite{JH98a}. In \cref{sec5}, we give concluding remarks.

\section{The GMTs and their inverses}\label{sec2}
M\"obius transformations, i.e., invertible rational mappings of degree $1$, have been widely used and extensively studied in complex analysis. In the theories of fields and projective spaces, a M\"obius transformation $P(x)=\frac{ax+b}{cx+d}$, where $a,b,c,d\in \FF{q}{n}$, $c\neq 0$ and $ad-bc\neq 0$, can be viewed as a mapping from $\FF{q}{n}\sqcup \{\infty\}$ to itself such that $P\bigl(-\frac{d}{c}\bigr)=\infty$ and $P(\infty)=\frac{a}{c}$. If $n=2$, we recall from \cref{sec1} that when applying the AGW-criterion over $\FF{q}{2}$, we usually employ those M\"obius transformations from $\PG(1,q)\cong \FP{q}\sqcup \{\infty\}$ to $\mu_{q+1}$ which are constructed as follows. Let $u\in \mu_{q+1}\setminus \FP{q}$, and define $\eta:\PG(1,q)\to \mu_{q+1}$, where
$$\eta(x)=\begin{cases}
1 & \text{if }x=\infty; \\
\frac{x+u^q}{x+u} & \text{if }x\in \FP{q}.
\end{cases}$$
Then $\eta$ is well-defined since $x+u\neq 0$ for all $x\in \FP{q}$ due to the assumption that $u\not\in \FP{q}$. It is also easy to see that $\eta$ is invertible and
$$\eta^{-1}(x)=\begin{cases}
\infty & \text{if }x=1; \\
\frac{-ux+u^q}{x-1} & \text{if }x\in \mu_{q+1}\setminus \{1\}.
\end{cases}$$
If, instead of identifying $\PG(1,q)$ with $\FP{q}\sqcup \{\infty\}$, we write it as $\{(x_0:1):x\in \FP{q}\}\sqcup \{(1:0)\}$ according to its original definition, and view $(1:0)$ as $\infty$, then for all $(x_0:x_1)\in \PG(1,q)$,
\begin{equation}\label{Eq005}
\eta(x_0:x_1)=(x_0+x_1u)^{q-1}=\frac{x_0+x_1u^q}{x_0+x_1u}.
\end{equation}
We also have that
\begin{equation}\label{Eq006}
\eta^{-1}(x)=\begin{cases}
(1:0) & \text{if }x=1; \\
\bigl(\frac{-ux+u^q}{x-1}:1\bigr) & \text{if }x\in \mu_{q+1}\setminus \{1\}.
\end{cases}
\end{equation}
These observations prompt us to introduce the notion of generalized M\"obius transformations (or GMTs for brevity), which are a class of bijections from $\PG(n-1,q)$ to $\mu_{\frac{q^n-1}{q-1}}$. 

\begin{dfn}\label{D21}
Let $W=\left\{\omega_j:0\leq j\leq n-1\right\}$ be a basis for $\FF{q}{n}$ over $\FP{q}$, and define $\psi_W:\PG(n-1,q)\to \mu_{\frac{q^n-1}{q-1}}$, where for each $(x_0:\dots:x_{n-1})\in \PG(n-1,q)$,
\begin{equation}\label{Eq007}
\psi_W(x_0:\dots:x_{n-1})
=\Biggl(\sum_{j=0}^{n-1}x_j\omega_j\Biggr)^{q-1}
=\frac{\displaystyle{\sum_{j=0}^{n-1}x_j\omega_j^q}}{\displaystyle{\sum_{j=0}^{n-1}x_j\omega_j}}.
\end{equation}
We call $\psi_W$ the generalized M\"obius transformation (GMT) associated with $W$.
\end{dfn}

A GMT $\psi_W$ is always well-defined. Indeed, $\psi_W(\PG(n-1,q))\subseteq \mu_{\frac{q^n-1}{q-1}}$ since
\begin{equation}\label{Eq008}
\mu_{\frac{q^n-1}{q-1}}=\left\{x^{q-1}:x\in \FF{q}{n}^\ast\right\}.
\end{equation}
Also, $\psi_W(\lambda x_0:\dots:\lambda x_{n-1})=\psi_W(x_0:\dots:x_{n-1})$ for each $(x_0:\dots:x_{n-1})\in \PG(n-1,q)$ and $\lambda\in \FP{q}^\ast$ because $\lambda^{q-1}=1$. Moreover, the denominator of the fractional form of $\psi_W$ in \cref{Eq007} is always non-zero since for each $(x_0:\dots:x_{n-1})\in \PG(n-1,q)$, at least one of $x_0,\dots,x_{n-1}$ is non-zero, meaning that $x_0\omega_0+\dots+x_{n-1}\omega_{n-1}\neq 0$ due to the $\FP{q}$-linearly independence of $W$. The following theorem shows that each GMT is a bijection.

\begin{thm}\label{T22}
For each basis $W$ for $\FF{q}{n}$ over $\FP{q}$, the mapping $\psi_W$ given by \cref{Eq007} is a bijection from $\PG(n-1,q)$ to $\mu_{\frac{q^n-1}{q-1}}$.
\end{thm}

\noindent \textit{Proof} Since $W$ is a basis for $\FF{q}{n}$ over $\FP{q}$, \cref{Eq007,Eq008} indicate that $\psi_W$ is a surjection. Since $\left|\PG(n-1,q)\right|=\Bigl|\mu_{\frac{q^n-1}{q-1}}\Bigr|$, $\psi_W$ is a bijection.\mbox{}\hfill $\square$\\
\\
When $n=2$, and $W=\{1,u\}$, where $u\in \mu_{q+1}\setminus \FP{q}$, the GMT $\psi_W$ defined by \cref{Eq007} agrees with the usual M\"obius transformation defined by \cref{Eq005}. When $n=3$ and $W=\left\{\omega,\omega^q,\omega^{q^2}\right\}$ for some $\omega\in \FF{q}{3}^\ast$,
$$\psi_W(x_0:x_1:x_2)=\frac{x_0\omega^q+x_1\omega^{q^2}+x_2\omega}{x_0\omega+x_1\omega^q+x_2\omega^{q^2}},$$
which agrees with the mapping in \cref{Eq002}, which was introduced in \cite{QL23a}.

In what follows, we determine, for each $n\geq 2$ and each basis $W$ for $\FF{q}{n}$ over $\FP{q}$, the inverse of the GMT $\psi_W$ defined by \cref{Eq007} using $W$, thereby explaining the known inverse formulas \cref{Eq006} (when $n=2$) and \cref{Eq003,Eq004} (when $n=3$). To achieve these requires a lemma about the dual bases.

\begin{lem}[{\cite[Theorem 8.7.]{Wan03}}]\label{L23}
Let $x\in \FF{q}{n}$. Let $W=\{\omega_j:0\leq j\leq n-1\}$ be a basis of $\FF{q}{n}$ over $\FP{q}$, and let $B=\{\beta_j:0\leq j\leq n-1\}$ be the dual basis of $W$. Then
\begin{equation}\label{Eq009}
x=\sum_{j=0}^{n-1}\Tr_q^{q^n}\left(\beta_jx\right)\omega_j.
\end{equation}
\end{lem}

To utilize \cref{Eq009}, we need to be able to find the dual basis $B$ of a given basis $W$. In general, this can be done by employing \textit{Moore matrices}. For an $n$-element subset $S=\left\{s_j:0\leq j\leq n-1\right\}$ of $\FF{q}{n}$, the Moore matrix associated with $S$ is
$$M_S
=\begin{bmatrix}
s_0     & s_0^q     & \dots  & s_0^{q^{n-1}} \\
s_1     & s_1^q     & \dots  & s_1^{q^{n-1}} \\
\vdots  & \vdots    & \ddots & \vdots \\
s_{n-1} & s_{n-1}^q & \dots  & s_{n-1}^{q^{n-1}} \\
\end{bmatrix}.$$

As a side note, Moore matrices find their applications in engineering and coding theory, especially when it comes to MRD codes. See for example \cite{BZ20a,CMPZ20a,CS19a,KL20a,KG05a,AN20a}. If $B$ and $W$ are the dual bases of each other and $I_n$ is the $n\times n$ identity matrix, then
$$M_WM_B^T
=\begin{bmatrix}
\omega_0     & \omega_0^q     & \dots  & \omega_0^{q^{n-1}} \\
\omega_1     & \omega_1^q     & \dots  & \omega_1^{q^{n-1}} \\
\vdots       & \vdots         & \ddots & \vdots \\
\omega_{n-1} & \omega_{n-1}^q & \dots  & \omega_{n-1}^{q^{n-1}} \\
\end{bmatrix}
\begin{bmatrix}
\beta_0           & \beta_1           & \dots  & \beta_{n-1} \\
\beta_0^q         & \beta_1^q         & \dots  & \beta_{n-1}^q \\
\vdots            & \vdots            & \ddots & \vdots \\
\beta_0^{q^{n-1}} & \beta_1^{q^{n-1}} & \dots  & \beta_{n-1}^{q^{n-1}} \\
\end{bmatrix}
=I_n.$$
In particular, $M_W$ is non-singular. Moreover, for $0\leq i\leq n-1$,
\begin{equation}\label{Eq010}
\beta_i=\det(M_W)^{-1}\cof_{M_W}(i,0),
\end{equation}
where $\cof_{M_W}(i,0)$ denotes the $(i,0)$-cofactor of $M_W$. The determinant of a Moore matrix has the following properties, which sometimes make determining dual bases via \cref{Eq010} less computationally expensive. 

\begin{prp}[\cite{BZ20a,DG96a}]\label{P24}
For any subset $S=\left\{s_j:0\leq j\leq n-1\right\}$ of $\FF{q}{n}$,
$$\det(M_S)=\prod_{j=0}^{n-1}\prod_{x_0,\dots,x_{j-1}\in \FP{q}}\left(\left(\sum_{i=0}^{j-1}x_is_i\right)+s_j\right).$$
In particular, $M_S$ is non-singular if and only if $S$ is a basis for $\FF{q}{n}$ over $\FP{q}$.
\end{prp}

\begin{prp}\label{P25}
If $W=\left\{\omega_j:0\leq j\leq n-1\right\}$ is a basis for $\FF{q}{n}$ over $\FP{q}$, then
\begin{equation}\label{Eq011}
\det(M_W)^{q-1}=(-1)^{n-1}.
\end{equation}
In particular, if $n$ is odd, then $\det(M_W)\in \FP{q}^\ast$.
\end{prp}

\noindent \textit{Proof}\quad Since $W$ is a basis for $\FF{q}{n}$ over $\FP{q}$,
$$\FF{q}{n}^\ast=\bigsqcup_{j=0}^{n-1}\left\{\lambda(x_0\omega_0+\dots+x_{j-1}\omega_{j-1}+\omega_j):\lambda\in \FP{q}^\ast\text{ and }x_0,\dots,x_{j-1}\in \FP{q}\right\}.$$
Hence, every element in $\FF{q}{n}^\ast$ is of the form $\lambda(x_0\omega_0+\dots+x_{j-1}\omega_{j-1}+\omega_j)$ for some unique $0\leq j\leq n-1,\lambda\in \FP{q}^\ast$ and $x_0,\dots,x_{j-1}\in \FP{q}$. Since the product of all non-zero elements in a finite field is $-1$, \cref{P24} implies that
\begin{align*}
\det(M_W)^{q-1}
& =1^n\det(M_W)^{q-1} \\
& =\biggl(-\prod_{\lambda\in \FP{q}^\ast}\lambda\biggr)^n\left(\prod_{j=0}^{n-1}\prod_{x_0,\dots,x_{j-1}\in \FP{q}}\left(\sum_{i=0}^{j-1}x_i\omega_i\right)+\omega_j\right)^{q-1} \\
& =(-1)^n\prod_{j=0}^{n-1}\prod_{\lambda\in \FP{q}^\ast}\prod_{x_0,\dots,x_{j-1}\in \FP{q}}\lambda\left(\left(\sum_{i=0}^{j-1}x_i\omega_i\right)+\omega_j\right) \\
& =(-1)^n\prod_{a\in \FF{q}{n}^\ast}a \\
& =(-1)^{n-1}.
\end{align*}
Thus, if $n$ is odd, then $\det(M_W)^{q-1}=1$, meaning that $\det(M_W)\in \FP{q}^\ast$. \mbox{}\hfill $\square$\\
\\
Let $\psi_W$ be the GMT associated with a basis $W$ for $\FF{q}{n}$ over $\FP{q}$. In view of the above results, we now give an inverse formula for $\psi_W$ which do not require us to actually find the dual basis $B$ of $W$.

\begin{thm}\label{T26}
Let $W=\{\omega_j:0\leq j\leq n-1\}$ be a basis for $\FF{q}{n}$ over $\FP{q}$. For $0\leq i\leq n-1$, let $c_i$ be the $(i,0)$-cofactor of $M_W$, the Moore matrix  associated with $W$. Then define
\begin{equation}\label{Eq012}
T_{W,i}(x)=\sum_{k=0}^{n-1}(-1)^{k(n-1)}c_i^{q^k}x^{\frac{q^k-1}{q-1}}.
\end{equation}
If $B=\{\beta_j:0\leq j\leq n-1\}$ is the dual basis of $W$ and $M_W$ is the Moore matrix associated with $W$, then for each $y\in \FF{q}{n}^\ast$,
\begin{equation}\label{Eq013}
\Tr_q^{q^n}(\beta_iy)=\det(M_W)^{-1}yT_{W,i}\left(y^{q-1}\right).
\end{equation}
In particular, $T_{W,i}\ (0\leq i\leq n-1)$ have no common roots in $\mu_{\frac{q^n-1}{q-1}}$, and for each $x\in \mu_{\frac{q^n-1}{q-1}}$,
\begin{equation}\label{Eq014}
\resizebox{0.8\hsize}{!}{%
$\psi_W^{-1}(x)
=\begin{cases}
(1:0:\dots:0) & \text{ if }x=\omega_0^{q-1}; \\
\displaystyle{\biggl(\frac{T_{W,0}(x)}{T_{W,1}(x)}:1:0:\dots:0\biggr)} & \text{ if }x=(x_0\omega_0+\omega_1)^{q-1}\text{ for some }x_0\in \FP{q};\\
\vdots &\ \vdots \\
\displaystyle{\biggl(\frac{T_{W,0}(x)}{T_{W,n-1}(x)}:\dots:\frac{T_{W,n-2}(x)}{T_{W,n-1}(x)}:1\biggr)} & \begin{array}{l}\text{if }x=(x_0\omega_0+\dots+x_{n-2}\omega_{n-2}+\omega_{n-1})^{q-1} \\ \text{for some }x_0,\dots,x_{n-2}\in \FP{q}.\end{array}
\end{cases}$%
}
\end{equation}
\end{thm}

\noindent \textit{Proof}\quad First, we verify \cref{Eq013}, using which we show that $T_{W,i}\ (0\leq i\leq n-1)$ have no common roots in $\mu_{\frac{q^n-1}{q-1}}$. For each $y\in \FF{q}{n}^\ast$,
\begin{align*}
\Tr_q^{q^n}(\beta_iy)\
& \stackrel{\mathmakebox[\widthof{=}]{(\ref{Eq010})}}{=}\ \sum_{k=0}^{n-1}\left(\det(M_W)^{-1}\cof_{M_W}(i,0)\right)^{q^k}y^{q^k} \\
& \stackrel{\mathmakebox[\widthof{=}]{(\ref{Eq011})}}{=}\ \det(M_W)^{-1}\sum_{k=0}^{n-1}(-1)^{k(n-1)}c_i^{q^k}y^{q^k} \\
& \stackrel{\mathmakebox[\widthof{=}]{(\ref{Eq012})}}{=}\ \det(M_W)^{-1}yT_{W,i}\left(y^{q-1}\right).
\end{align*}
Hence, \cref{Eq013} holds. By \cref{Eq009}, $\Tr_q^{q^n}(\beta_iy)\ (0\leq i\leq n-1)$ have no common roots in $\FF{q}{n}^\ast$. Therefore, $T_{W,i}\ (0\leq i\leq n-1)$ have no common roots in $\mu_{\frac{q^n-1}{q-1}}$.

Next, we prove \cref{Eq014}. Assume that $y\in \FF{q}{n}^\ast$ and $x=y^{q-1}$. Let $0\leq j\leq n-1$. Owing to \cref{Eq009,Eq013}, the following are equivalent. 
\begin{enumerate}
    \item $j=\max\left\{0\leq i\leq n-1:T_{W,i}(x)\neq 0\right\}$;
    \item $j=\max\left\{0\leq i\leq n-1:\Tr_q^{q^n}(\beta_iy)\neq 0\right\}$;
    \item $y=\lambda (x_0\omega_0+\dots+x_{j-1}\omega_{j-1}+\omega_j)$ for some $x_0,\dots,x_{j-1}\in \FP{q}$ and $\lambda\in \FP{q}^\ast$;
    \item $x=(x_0\omega_0+\dots+x_{j-1}\omega_{j-1}+\omega_j)^{q-1}$ for some $x_0,\dots,x_{j-1}\in \FP{q}$.
\end{enumerate}
Hence, if $j=\max\left\{0\leq i\leq n-1:\Tr_q^{q^n}(\beta_iy)\neq 0\right\}$, then
\begin{align*}
\psi_W^{-1}(x)
& =\psi_W^{-1}\left(y^{q-1}\right) \\
& \stackrel{\mathmakebox[\widthof{=}]{(\ref{Eq009})}}{=}\left(\Tr_q^{q^n}(\beta_0y):\dots:\Tr_q^{q^n}(\beta_{n-1}y)\right) \\
& =\Biggl(\frac{\Tr_q^{q^n}(\beta_0y)}{\Tr_q^{q^n}(\beta_jy)}:\dots:\frac{\Tr_q^{q^n}(\beta_{j-1}y)}{\Tr_q^{q^n}(\beta_jy)}:1:0:\dots:0\Biggr) \\
& \stackrel{\mathmakebox[\widthof{=}]{(\ref{Eq013})}}{=}\biggl(\frac{T_{W,0}(x)}{T_{W,j}(x)}:\dots:\frac{T_{W,j-1}(x)}{T_{W,j}(x)}:1:0:\dots:0\biggr),
\end{align*}
which proves \cref{Eq014}.\mbox{}\hfill $\square$\\
\\
For each $T_{W,i}$, all coefficients can be determined by $c_i=T_{W,i}(0)$. Since $c_i$ is the $(i,0)$-cofactor of $M_W$, we obtain from the Leibniz formula for determinants \cite{ANP16a} that
$$c_i=(-1)^i\left(\sum_{\sigma\in S_n}\sgn(\sigma)\prod_{k=0}^{i-1}\omega_k^{\sigma(k)}\prod_{k=i}^{n-2}\omega_{k+1}^{\sigma(k)}\right)^q,$$
where $S_n$ is the symmetric group consisting of the $n!$ permutations of $\{0,1,\dots,n-1\}$.\\
\\
When $n=2$ and $W=\{1,u\}$, where $u\in \mu_{q+1}\setminus \FP{q}$, it is easy to check that $T_{W,0}(x)=-ux+u^q$ and $T_{W,1}(x)=x-1$. Thus, the inverse formula \cref{Eq006} for a usual M\"obius transformation is a special case of \cref{Eq014}.

When $W$ is \textit{normal}, the process of finding $\psi_W^{-1}$ via \cref{T26} can be further simplified. Indeed, we recall that a basis $W$ for $\FF{q}{n}$ over $\FP{q}$ is normal when $W=\bigl\{\omega^{q^j}:0\leq j\leq n-1\bigr\}$ for some $\omega\in \FF{q}{n}^\ast$. In that case, \cite[Theorem 2.2.7.]{SG93a} shows that the dual basis $B$ of $W$ is also normal so that $B=\left\{\beta^{q^i}:0\leq i\leq n-1\right\}$ for some $\beta\in \FF{q}{n}^\ast$. Hence, from \cref{Eq010,Eq011}, we deduce that for $0\leq i\leq n-1$,
\begin{align*}
c_i
& =\cof_{M_W}(i,0) \\
& =\det(M_W)\beta^{q^i} \\
& =\det(M_W)\left(\det(M_W)^{-1}\cof_{M_W}(0,0)\right)^{q^i} \\
& =\det(M_W)^{1-q^i}\cof_{M_W}(0,0)^{q^i} \\
& =(-1)^{(n-1)i}c_0^{q^i}.
\end{align*}
So finding $\psi_W^{-1}$ boils down to determining $c_0$, which is much less computationally expensive than computing all of $c_i\ (0\leq i\leq n-1)$. For example, let $n=3$ and $W=\left\{\omega,\omega^q,\omega^{q^2}\right\}$ a normal basis for $\FF{q}{3}$ over $\FP{q}$. A direct calculation shows that $\det(M_W)=3\N_q^{q^3}(\omega)-\Tr_q^{q^3}\left(\omega^3\right)\in \FP{q}^\ast$, and the dual basis of $W$ is $B=\left\{\beta,\beta^q,\beta^{q^2}\right\}$, where $\beta=\det(M_W)^{-1}\bigl(\omega^{q^2+q}-\omega^2\bigr)$. Hence, letting $\alpha=\omega^{q+1}-\omega^{2q^2}$ yields that $T_{W,i}(x)=\alpha^{q^{i+1}}+\alpha^{q^{i+2}}x+\alpha^{q^{i+3}}x^{q+1}\ (i=0,1,2)$, which agree with the polynomials given by \cref{Eq003}. In this case, $\psi_W^{-1}$, which is given by \cref{Eq004}, is a special case of \cref{Eq014}.

We also remark that if $W$ is a \textit{polynomial basis}, i.e., $W=\bigl\{\omega^j:0\leq j\leq n-1\bigr\}$ for some $\omega\in \FF{q}{n}^\ast$, an easy way to compute the dual basis of $W$ is given by \cite[Theorem 8.6.]{Wan03}.

\section{Applications}\label{sec3}
So far, we have generalized the M\"obius transformations and determined the inverse of the GMTs. Our strategy for involving the GMTs in constructing permutation polynomials of the form $f(x)=x^rh\bigl(x^{q-1}\bigr)$ over $\FF{q}{n}$ is divide-and-conquer. We construct $h$ and $f$ piecewise based on the following partitions.

\begin{thm}\label{T31}
Let $W=\{\omega_j:0\leq j\leq n-1\}$ be any basis for $\FF{q}{n}$ over $\FP{q}$, and let $B=\{\beta_j:0\leq j\leq n-1\}$ be the dual basis of $W$. Let $T_{W,i}\ (0\leq i\leq n-1)$ be as defined by \cref{Eq012}. For $0\leq j\leq n-1$, define
\begin{align*}
S_{W,j} & =\bigl\{x\in \FF{q}{n}^\ast:j=\max\bigl\{0\leq i\leq n-1:\Tr_q^{q^n}(\beta_ix)\neq 0\bigr\}\bigr\}; \\
Z_{W,j} & =\bigl\{x\in \mu_{\frac{q^n-1}{q-1}}:j=\max\bigl\{0\leq i\leq n-1:T_{W,j}(x)\neq 0\bigr\}\bigr\}; \\
\overline{C_{q,n,j}} & =\biggl\{\biggl(\frac{Tr_q^{q^n}(\beta_0x)}{Tr_q^{q^n}(\beta_jx)}:\dots:\frac{Tr_q^{q^n}(\beta_{j-1}x)}{Tr_q^{q^n}(\beta_jx)}:1:0:\dots:0\biggr):x\in S_{W,j}\biggr\}.
\end{align*}
Then we have the following partitions.
\begin{align*}
\FF{q}{n}^\ast=\bigsqcup_{j=0}^{n-1}S_{W,j};\quad
\mu_{\frac{q^n-1}{q-1}}=\bigsqcup_{j=0}^{n-1}Z_{W,j};\quad
\PG(n-1,q)=\bigsqcup_{j=0}^{n-1}\overline{C_{q,n,j}}.
\end{align*}
Moreover, for $0\leq j\leq n-1$,
\begin{equation}\label{Eq015}
Z_{W,j}=\psi_W\bigl(\overline{C_{q,n,j}}\bigr)=S_{W,j}^{q-1}.
\end{equation}  
\end{thm}

\noindent \textit{Proof}\quad Since $\Tr_q^{q^n}(\beta_ix)\ (0\leq i\leq n-1)$ have no common roots in $\FF{q}{n}^\ast$, the disjoint union of $S_{W,j}\ (0\leq j\leq n-1)$ is $\FF{q}{n}^\ast$. According to items $(2)$ and $(4)$ of the four equivalence relations in the last part of the proof of \cref{T26}, the disjoint union of $Z_{W,j}\ (0\leq j\leq n-1$ is $\mu_{\frac{q^n-1}{q-1}}$. To show that the disjoint union of $\overline{C_{q,n,j}}\ (0\leq j\leq n-1)$ is $\PG(n-1,q)$, we notice that for any given elements $x_0,\dots,x_{j-1}\in \FP{q}$, if $x=x_0\omega_0+\dots+x_{j-1}\omega_{j-1}+\omega_j$, then $\Tr_q^{q^n}(\beta_jx)=1$, and $\Tr_q^{q^n}(\beta_ix)=x_i$ if $0\leq i\leq j-1$ and $0$ otherwise. Hence, for $0\leq j\leq n-1$,
\begin{equation}\label{Eq016}
\overline{C_{q,n,j}}=\{(x_0:\dots:x_{j-1}:1:0:\dots:0):x_0,\dots,x_{j-1}\in \FP{q}\},
\end{equation}
Thus, their disjoint union is $\PG(n-1,q)$. \cref{Eq015} follows directly from \cref{Eq007,Eq016} and the four equivalence relations in the last part of the proof of \cref{T26}.\mbox{}\hfill $\square$

\begin{figure}[hbt!]
\begin{tikzcd}[row sep=4ex, column sep=14ex]
{\displaystyle{\bigsqcup_{j=0}^{n-1}S_{W,j}}} \arrow{r}{f(x)=x^rh\left(x^{q-1}\right)} \arrow{d}{x^{q-1}} & {\displaystyle{\bigsqcup_{j=0}^{n-1}S_{W,j}}} \arrow{d}{x^{q-1}} \\
{\displaystyle{\bigsqcup_{j=0}^{n-1}Z_{W,j}}} \arrow{r}{g(x)=x^rh(x)^{q-1}} \arrow{d}{\psi_W^{-1}} & {\displaystyle{\bigsqcup_{j=0}^{n-1}Z_{W,j}}} \\
{\displaystyle{\bigsqcup_{j=0}^{n-1}\overline{C_{q,n,j}}}} \arrow{r}{\overline{g}} & {\displaystyle{\bigsqcup_{j=0}^{n-1}\overline{C_{q,n,j}}}} \arrow{u}[right]{\psi_Y} 
\end{tikzcd}{\caption{The multiplicative AGW criterion with explicit partitions}\label{Fig3}}
\end{figure}

Let $W$ and $Y=\{\gamma_j:0\leq j\leq n-1\}$ be any two bases for $\FF{q}{n}$ over $\FP{q}$, and let $B=\left\{\beta_j:0\leq j\leq n-1\right\}$ be the dual basis of $W$. Define $\psi_W,\psi_Y$ by \cref{Eq007}, and let $\overline{g}$ be any bijection over $\PG(n-1,q)$. Clearly, $\psi_Y\circ \overline{g}\circ \psi_W^{-1}$ is a bijection over $\mu_{\frac{q^n-1}{q-1}}$. In order to construct permutation polynomials of the form $f(x)=x^rh\bigl(x^{q-1}\bigr)$ over $\FF{q}{n}$, where $\gcd(r,q-1)=1$, we want to find examples of bijections $\overline{g}$ over $\PG(n-1,q)$ and $h\in \FF{q}{n}[x]$ such that for all $x\in \mu_{\frac{q^n-1}{q-1}}$,
\begin{equation}\label{Eq017}
x^rh(x)^{q-1}=\psi_Y\circ \overline{g}\circ \psi_W^{-1}(x).
\end{equation}
As \cref{Fig3} shows, it suffices to ensure that \cref{Eq017} holds on each of the sets $Z_{W,j}$ because their disjoin union is $\mu_{\frac{q^n-1}{q-1}}$ by \cref{T31}. We observe that each mapping $\overline{g}$ from $\PG(n-1,q)$ to itself can be written as $\overline{g}=(g_0:\dots:g_{n-1})$ for some $g_0,\dots,g_{n-1}\in \FP{q}[x_0,\dots,x_{n-1}]$. Hence, for each $(x_0,\dots,x_{n-1})\in \FP{q}^n\setminus \bigl\{\vec{0}\bigr\}$,
$$\overline{g}(x_0:\dots:x_{n-1})=(g_0(x_0,\dots,x_{n-1}):\dots:g_{n-1}(x_0,\dots,x_{n-1})).$$
We therefore want bijections $\overline{g}$ over $\PG(n-1,q)$ and $h\in \FF{q}{n}[x]$ such that for each $j$,
\begin{equation}\label{Eq018}
\begin{aligned}
x^rh(x)^{q-1}\
& \stackrel{\mathmakebox[\widthof{=}]{(\ref{Eq017})}}{=}\ \psi_Y\circ \overline{g}\circ \psi_W^{-1}(x) \\
& \stackrel{\mathmakebox[\widthof{=}]{(\ref{Eq007},\ref{Eq014})}}{=}\ \biggl(\sum_{k=0}^{n-1}g_k\biggl(\frac{T_{W,0}(x)}{T_{W,j}(x)},\dots,\frac{T_{W,j-1}(x)}{T_{W,j}(x)},1,0,\dots,0\biggr)\gamma_k\biggr)^{q-1}    
\end{aligned}
\end{equation}
for all $x\in Z_{W,j}$. For each such pair $\overline{g},h$, we know that $f(x)=x^rh\bigl(x^{q-1}\bigr)$ permutes $\FF{q}{n}$.\\
\\
We also note that in \cref{T31}, we use $\Tr_q^{q^n}(\beta_ix)\ (0\leq i\leq n-1)$ to partition $\FF{q}{n}^\ast$. But according to the definition of dual basis, for $0\leq j\leq n-1$, the set $S_{W,j}$ defined in \cref{T31} can be rewritten as
\begin{equation}\label{Eq019}
S_{W,j}=\{\lambda(x_0\omega_0+\dots x_{j-1}\omega_{j-1}+\omega_j):\lambda\in \FP{q}^\ast\text{ and }x_0,\dots,x_{j-1}\in \FP{q}\}.
\end{equation}
Hence, we may use the two definitions interchangeably, depending on which one is more convenient. Based on \cref{Eq018}, the following theorem details a method of constructing permutations over $\FF{q}{n}$ using the GMTs and bijections over $\PG(n-1,q)$.

\begin{thm}\label{T32}
Let $W$ and $Y=\left\{\gamma_j:0\leq j\leq n-1\right\}$ be any two bases for $\FF{q}{n}$ over $\FP{q}$. Let $B=\left\{\beta_j:0\leq j\leq n-1\right\}$ be the dual basis of $W$. Let $\widetilde{g}=(g_0,\dots,g_{n-1})$ be any vector-valued polynomial such that $\overline{g}=(g_0:\dots:g_{n-1})$ is a bijection over $\PG(n-1,q)$. Let $r\in \NN$ be coprime to $q-1$. For $0\leq j\leq n-1$ and each $x\in \FF{q}{n}^\ast$ such that
$$j=\max\bigl\{0\leq i\leq n-1:\Tr_q^{q^n}(\beta_ix)\neq 0\bigr\},$$
we define
\begin{equation}\label{Eq020}
\begin{aligned}
f_j(x)
& =\Tr_q^{q^n}(\beta_jx)^r\sum_{k=0}^{n-1}g_k\Biggl(\frac{\Tr_q^{q^n}(\beta_0x)}{\Tr_q^{q^n}(\beta_jx)},\dots,\frac{\Tr_q^{q^n}(\beta_{j-1}x)}{\Tr_q^{q^n}(\beta_jx)},1,0,\dots,0\Biggr)\gamma_k. \\
& =\Tr_q^{q^n}(\beta_jx)^r\widetilde{g}\Biggl(\frac{\Tr_q^{q^n}(\beta_0x)}{\Tr_q^{q^n}(\beta_jx)},\dots,\frac{\Tr_q^{q^n}(\beta_{j-1}x)}{\Tr_q^{q^n}(\beta_jx)},1,0,\dots,0\Biggr)\cdot (\gamma_0,\dots,\gamma_{n-1}),
\end{aligned}
\end{equation}
where $\cdot$ denotes the dot product of two vectors in $\FP{q}^n$. Define $f:\FF{q}{n}\to \FF{q}{n}$, where
\begin{equation}\label{Eq021}
f(x)=\begin{cases}
0 & \text{if }x=0; \\
f_j(x) & \text{if }j=\max\bigl\{0\leq i\leq n-1:\Tr_q^{q^n}(\beta_ix)\neq 0\bigr\}\quad (0\leq j\leq n-1),
\end{cases}
\end{equation}
Then $f$ is an $r$-homogeneous rational mapping which permutes $\FF{q}{n}$. In particular, if $r\geq \max\{\deg(g_k):0\leq k\leq n-1\}$, then $f$ is an $r$-homogeneous permutation polynomial of index $\frac{q^n-1}{q-1}$ over $\FF{q}{n}$.
\end{thm}

\noindent \textit{Proof}\quad It follows trivially from \cref{Eq020} and the $\FP{q}$-linearity of $\Tr_q^{q^n}$ that each $f_j$ is $r$-homogeneous. Hence, $f$ is $r$-homogeneous due to \cref{Eq021}. Now we show that $f$ permutes $\FF{q}{n}$. Let $T_{W,i}\ (0\leq i\leq n-1)$ be as defined by \cref{Eq012}. Define $\phi_W:\mu_{\frac{q^n-1}{q-1}}\to \FP{q}^n\setminus \bigl\{\vec{0}\bigr\}$, where for $0\leq j\leq n-1$ and each $x\in Z_{W,j}$ (which is defined in \cref{T31}),
\begin{equation}\label{Eq022}
\phi_W(x)=\biggl(\frac{T_{W,0}(x)}{T_{W,j}(x)},\dots,\frac{T_{W,j-1}(x)}{T_{W,j}(x)},1,0,\dots,0\biggr).
\end{equation}
Let $M_W$ be the Moore matrix associated with $W$. For $0\leq j\leq n-1$, let
\begin{equation}\label{Eq023}
h_j(x)=\det(M_W)^{-r}T_{W,j}(x)^r\sum_{k=0}^{n-1}g_k\left(\phi_W(x)\right)\gamma_k.
\end{equation}
Let $\widetilde{h}\in \FF{q}{n}[x]$ be an arbitrary mapping, and define
\begin{equation}\label{Eq024}
h(x)=
\begin{cases}
\widetilde{h}(x) & \text{if }x\in \FF{q}{n}\setminus \mu_{\frac{q^n-1}{q-1}}; \\  
h_j(x) & \text{if }x\in Z_{W,j}\quad (0\leq j\leq n-1).
\end{cases}
\end{equation}
First, we show that $x^rh(x)^{q-1}$ permutes $\mu_{\frac{q^n-1}{q-1}}$. Then we prove that the mapping $f(x)$ defined by \cref{Eq020,Eq021} in fact equals $x^rh\left(x^{q-1}\right)$, hence permutes $\FF{q}{n}$. To verify that $x^rh(x)^{q-1}$ permutes $\mu_{\frac{q^n-1}{q-1}}$, it suffices to prove \cref{Eq017}. Assume that $x\in \mu_{\frac{q^n-1}{q-1}}$. By \cref{T31}, there exists a unique $0\leq j\leq n-1$ such that $x\in Z_{W,j}$. If $x=y^{q-1}$, where $y\in \FF{q}{n}^\ast$, then
\begin{equation}\label{Eq025}
\begin{aligned}
T_{W,j}(x)^q\
& \stackrel{\mathmakebox[\widthof{=}]{(\ref{Eq013})}}{=}\ \left(\det(M_W)y^{-1}\Tr_q^{q^{n}}(\beta_jy)\right)^q \\
& \stackrel{\mathmakebox[\widthof{=}]{(\ref{Eq011})}}{=}\ (-1)^{n-1}y^{-(q-1)}\left(\det(M_W)y^{-1}\Tr_q^{q^{n}}(\beta_jy)\right) \\
& \stackrel{\mathmakebox[\widthof{=}]{(\ref{Eq012})}}{=}\ (-1)^{n-1}x^{-1}T_{W,j}(x).   
\end{aligned}
\end{equation}
We therefore deduce that
\begin{align*}
\psi_Y\circ \overline{g}\circ \psi_W^{-1}(x)\ \
& \stackrel{\mathmakebox[\widthof{=}]{(\ref{Eq014},\ref{Eq022})}}{=}\ \ \Biggl(\sum_{k=0}^{n-1}g_k\left(\phi_W(x)\right)\gamma_k\Biggr)^{q-1} \\
& \stackrel{\mathmakebox[\widthof{=}]{(\ref{Eq023})}}{=}\ \ \bigl(\det(M_W)^rT_{W,j}(x)^{-r}h_j(x)\bigr)^{q-1} \\
& \stackrel{\mathmakebox[\widthof{=}]{(\ref{Eq011},\ref{Eq024})}}{=}\ \ (-1)^{r(n-1)}T_{W,j}(x)^{(q-1)(-r)}h(x)^{q-1} \\
& \stackrel{\mathmakebox[\widthof{=}]{(\ref{Eq025})}}{=}\ \ x^rh(x)^{q-1}.
\end{align*}

So $x^rh(x)^{q-1}$ and $x^rh\left(x^{q-1}\right)$ permute $\mu_{\frac{q^n-1}{q-1}}$ and $\FF{q}{n}$, respectively. They are both rational since $T_{W,i},g_k\ (0\leq i,k\leq n-1)$ are polynomials and all components of $\phi_W$ are rational. If $r\geq \max\{\deg(g_k):0\leq k\leq n-1\}$, then for $0\leq j\leq n-1$ and each $x\in Z_{W,j}$, it is easy to see that $T_{W,j}(x)^rg_k\left(\phi_{W,n}(x)\right)\ (0\leq k\leq n-1)$ are polynomials. So $h_j\ (0\leq j\leq n-1)$ and $h$ are polynomials. Consequently, $x^rh(x)^{q-1}$ and $x^rh\left(x^{q-1}\right)$ are both polynomials in this case. Finally, we show that the mapping $f$ defined by \cref{Eq020,Eq021} is a permutation polynomial of index $\frac{q^n-1}{q-1}$ over $\FF{q}{n}$ by verifying that
$$f(x)=x^rh\left(x^{q-1}\right).$$
This is clearly true when $x=0$. If $x\in \FF{q}{n}^\ast$, then there exists a unique $0\leq j\leq n-1$ such that $x\in S_{W,j}$. That is, $j=\max\bigl\{0\leq i\leq n-1:\Tr_q^{q^n}(\beta_ix)\neq 0\bigr\}$. Consequently,
\begin{align*}
f(x)\ \
\stackrel{\mathmakebox[\widthof{=}]{(\ref{Eq020},\ref{Eq021})}}{=}\ &\ \ Tr_q^{q^n}(\beta_jx)^r\sum_{k=0}^{n-1}g_k\biggl(\frac{Tr_q^{q^n}(\beta_0x)}{Tr_q^{q^n}(\beta_jx)},\dots,\frac{Tr_q^{q^n}(\beta_{j-1}x)}{Tr_q^{q^n}(\beta_jx)},1,0,\dots,0\biggr)\gamma_k. \\
\stackrel{\mathmakebox[\widthof{=}]{(\ref{Eq013},\ref{Eq022})}}{=}\ &\ \ \bigl(\det(M_W)^{-1}xT_{W,j}\bigl(x^{q-1}\bigr)\bigr)^r\sum_{k=0}^{n-1}g_k\bigl(\phi_W\bigl(x^{q-1}\bigr)\bigr)\gamma_k \\
\stackrel{\mathmakebox[\widthof{=}]{(\ref{Eq023},\ref{Eq024})}}{=}\ &\ \ x^rh\bigl(x^{q-1}\bigr).
\end{align*}
Thus, $f$ is a permutation polynomial of index $\frac{q^n-1}{q-1}$ over $\FF{q}{n}$.\mbox{}\hfill $\square$\\
\\
Before demonstrating any applications of \cref{T32}, there are a couple of important yet unanswered questions. That is, when is a mapping $\overline{g}$ over $\PG(n-1,q)$ well-defined, and when is it a bijection? To answer these questions, we need the following results.

\begin{lem}\label{L33}
Let $\overline{g}=(g_0:\dots:g_{n-1})$ be an arbitrary mapping from $\PG(n-1,q)$ to itself, where $g_0,\dots,g_{n-1}:\FP{q}^n\to \FP{q}$. Let $\widetilde{g}=(g_0,\dots,g_{n-1})$. Then $\overline{g}$ is well-defined if and only if
\begin{enumerate}
    \item $g_0,\dots,g_{n-1}$ have no common roots in $\FP{q}^n\setminus \bigl\{\vec{0}\bigr\}$ and \label{L33I1}
    \item there exists a mapping $D:\FP{q}^\ast\times \bigl(\FP{q}^n\setminus \bigl\{\vec{0}\bigr\}\bigr)\to \FP{q}^\ast$ such that
    $$\widetilde{g}\bigl(\lambda \vec{x}\bigr)=D\bigl(\lambda,\vec{x}\bigr)\widetilde{g}\bigl(\vec{x}\bigr)$$
    for all $\lambda\in \FP{q}^\ast$ and $\vec{x}\in \FP{q}^n\setminus \bigl\{\vec{0}\bigr\}$. \label{L33I2}
\end{enumerate}
\end{lem}

We skip the proof of \cref{L33}, which is a routine check. However, \cref{L33} gives rise to a criterion for determining when a mapping $\overline{g}=(g_0:\dots:g_{n-1})$ permutes $\PG(n-1,q)$ based on the images of the underlying sets of all equivalence classes in $\PG(n-1,q)$ under the mapping $\widetilde{g}=(g_0,\dots,g_{n-1})$.

\begin{prp}\label{P34}
Let $\overline{E_i},E_i\ \left(0\leq i\leq \frac{q^n-1}{q-1}-1\right)$ be all equivalence classes in $\PG(n-1,q)$ and their underlying sets, respectively. That is, each $E_i$ is the set of all representatives for $\overline{E_i}$ in $\FP{q}^n\setminus \bigl\{\vec{0}\bigr\}$. Let $\overline{g}=(g_0:\dots:g_{n-1})$ be any well-defined mapping from $\PG(n-1,q)$ to itself, where $g_0,\dots,g_{n-1}:\FP{q}^n\to \FP{q}$. Let $\widetilde{g}=(g_0,\dots,g_{n-1})$. Then $\overline{g}$ is a bijection over $\PG(n-1,q)$ if and only if there exists a permutation $\sigma$ in the symmetric group $S_{\frac{q^n-1}{q-1}}$ such that $\widetilde{g}(E_i)\subseteq E_{\sigma(i)}$ for $0\leq i\leq \frac{q^n-1}{q-1}-1$.
\end{prp}

\noindent \textit{Proof}\quad We show that for $0\leq i_1,i_2\leq \frac{q^n-1}{q-1}-1$, $\widetilde{g}\left(E_{i_1}\right)\subseteq E_{i_2}$ if and only if $\overline{g}\left(\overline{E_{i_1}}\right)=\overline{E_{i_2}}$. Assume that $E_{i_1}=\left\{\delta \vec{x}:\delta\in \FP{q}^\ast\right\}$ and $E_{i_2}=\left\{\delta \vec{y}:\delta\in \FP{q}^\ast\right\}$, where $\vec{x}=(x_0,\dots,x_{n-1})\in E_{i_1}$ and $\vec{y}=(y_0,\dots,y_{n-1})\in E_{i_2}$. Since $\overline{g}$ is well-defined, \cref{L33I2} of \cref{L33} implies that
\begin{align*}
\widetilde{g}\left(E_{i_1}\right)\subseteq E_{i_2}
&\ \Leftrightarrow\ \widetilde{g}\left(\vec{x}\right)=\left(g_0\left(\vec{x}\right),\dots,g_{n-1}\left(\vec{x}\right)\right)=\delta \vec{y}\text{ for some }\delta\in \FP{q}^\ast \\
&\ \Leftrightarrow\ \overline{g}(x_0:\dots:x_{n-1})=\left(g_0\left(\vec{x}\right):\dots:g_{n-1}\left(\vec{x}\right)\right)=(y_0:\dots:y_{n-1}) \\
&\ \Leftrightarrow\ \overline{g}\left(\overline{E_{i_1}}\right)=\overline{E_{i_2}}.
\end{align*}
Consequently, $\overline{g}$ is a bijection over $\PG(n-1,q)$ if and only if $\widetilde{g}(E_i)\subseteq E_{\sigma(i)}$ for $0\leq i\leq \frac{q^n-1}{q-1}-1$ for some $\sigma\in S_{\frac{q^n-1}{q-1}}$.\mbox{}\hfill $\square$\\
\\
The major takeaway from \cref{P34} is that for $\overline{g}=(g_0:\dots:g_{n-1})$ to be a bijection over $\PG(n-1,q)$, $\widetilde{g}=(g_0,\dots,g_{n-1})$ need not be a bijection over $\FP{q}^n$. Instead, it suffices to make sure that $\widetilde{g}$ maps the underlying sets of equivalence classes in $\PG(n-1,q)$ to subsets of pairwise distinct underlying sets of equivalence classes in $\PG(n-1,q)$. In the next few results, we present some concrete examples of homogeneous permutation polynomials of index $\frac{q^n-1}{q-1}$ over $\FF{q}{n}$ which are constructed from bijections over $\PG(n-1,q)$.

\begin{prp}\label{P35}
Let $H_0(x)=x$, and let $H_j\ (1\leq j\leq n-1)$ be any polynomials over $\FF{q}{n}$ which permute $\FP{q}$. For $0\leq k\leq n-1$, define $g_k:\FP{q}^n\to \FP{q}$, where
$$g_k(x_0,\dots,x_{n-1})=\sum_{j=k}^{n-1}H_j\Bigl(x_kx_j^{q-2}\Bigr)x_j^{q-1}\prod_{i=j+1}^{n-1}\Bigl(1-x_i^{q-1}\Bigr).$$
Then $\overline{g}=(g_0:\dots:d_{n-1})$ permutes $\PG(n-1,q)$.
\end{prp}
\noindent \textit{Proof}\quad For $0\leq j,k\leq n-1$, it is easy to check that
$$g_k(x_0,\dots,x_{n-1})=
\begin{cases}
H_j\bigl(\frac{x_k}{x_j}\bigr) & \text{if }0\leq k\leq j,x_j\neq 0\text{ and }x_i=0\text{ for all }i\geq j+1; \\
0 & \text{otherwise}.
\end{cases}$$
Let $\widetilde{g}=(g_0,\dots,g_{n-1})$. Then $\widetilde{g}\bigl(\vec{0}\bigr)=\vec{0}$, and for all $\lambda\in \FP{q}^\ast$,
$$\lambda(1,0,\dots,0)\stackrel{\widetilde{g}}{\mapsto}(1,0,\dots,0).$$
Moreover, if $1\leq j\leq n-1,\lambda\in \FP{q}^\ast$ and $x_0,\dots,x_{j-1}\in \FP{q}$, then
$$\lambda (x_0,\dots,x_{j-1},1,0,\dots,0)\stackrel{\widetilde{g}}{\mapsto}(H_j(x_0),\dots,H_j(x_{j-1}),1,0,\dots,0).$$
Since $\vec{0}$ is the only root of $\widetilde{g}$ in $\FP{q}^n$, and $\widetilde{g}\bigl(\lambda\vec{x}\bigr)=\widetilde{g}\bigl(\vec{x}\bigr)$ for $0\leq k\leq n-1$, $\overline{g}=(g_0:\dots:g_{n-1})$ satisfies both conditions in \cref{L33}, hence is well-defined over $\PG(n-1,q)$. Also, we have that for $0\leq j\leq n-1$ and $x_0,\dots,x_{j-1}\in \FP{q}$,
\begin{equation}\label{Eq026}
\overline{C_{q,n,j}}\stackrel{(\ref{Eq016})}{=}(x_0:\dots:x_{j-1}:1:0:\dots:0)\stackrel{\overline{g}}{\mapsto}(H_j(x_0):\dots:H_j(x_{j-1}):1:0:\dots:0).
\end{equation}
Since $H_0,\dots,H_{n-1}$ permute $\FP{q}$, $\overline{g}$ permutes each of the sets $\overline{C_{q,n,j}}\ (0\leq j\leq n-1)$, the union of which is $\PG(n-1,q)$. Consequently, $\overline{g}$ permutes $\PG(n-1,q)$.\mbox{}\hfill $\square$\\
\\
We deduce from \cref{Eq026} that for $0\leq j\leq n-1$ and each $x\in S_{W,j}$,
\begin{align*}
g_k\biggl(\frac{\Tr_q^{q^n}(\beta_0x)}{\Tr_q^{q^n}(\beta_jx)},\dots,\frac{\Tr_q^{q^n}(\beta_{j-1}x)}{\Tr_q^{q^n}(\beta_jx)},1,0,\dots,0\biggr)
& =\begin{cases}
H_j\biggl(\frac{\Tr_q^{q^n}(\beta_kx)}{\Tr_q^{q^n}(\beta_jx)}\biggr) & \text{if }0\leq k\leq j-1; \\
1 & \text{if }k=j; \\
0 & \text{otherwise},
\end{cases}
\end{align*}
Plugging the above relations in \cref{Eq020} in \cref{T32} yields the following classes of $r$-homogeneous permutation polynomials of index $\frac{q^n-1}{q-1}$ over $\FF{q}{n}$.

\begin{thm}\label{T36}
Let $W$ and $Y=\left\{\gamma_j:0\leq j\leq n-1\right\}$ be any two bases for $\FF{q}{n}$ over $\FP{q}$. Let $B=\left\{\beta_j:0\leq j\leq n-1\right\}$ be the dual basis of $W$. Let $H_0(x)=x$, and let $H_j\ (1\leq j\leq n-1)$ be any polynomials over $\FF{q}{n}$ which permute $\FP{q}$. Let $r\in \NN$ be coprime to $q-1$, and assume that $r\geq \deg(H_j)\ (1\leq j\leq n-1)$. For $0\leq j\leq n-1$ and each $x\in S_{W,j}$, let
\begin{equation}\label{Eq027}
F_j(x)=\Tr_q^{q^n}(\beta_jx)^r\biggl(\sum_{k=0}^{j-1}H_j\biggl(\frac{\Tr_q^{q^n}(\beta_kx)}{\Tr_q^{q^n}(\beta_jx)}\biggr)\gamma_k+\gamma_j\biggr).
\end{equation}
Then the following mapping is an $r$-homogeneous permutation polynomial over $\FF{q}{n}$.
$$F(x)=\begin{cases}
0 & \text{if }x=0; \\
F_j(x) & \text{if }j=\max\bigl\{0\leq i\leq n-1:\Tr_q^{q^n}(\beta_ix)\neq 0\bigr\}\quad (0\leq j\leq n-1).
\end{cases}$$
\end{thm}

The case when $n=2$ in \cref{T36} is particularly interesting since in that case, we can construct homogeneous permutation polynomials of index $q+1$ over $\FF{q}{2}$ without explicitly picking two bases for $\FF{q}{2}$ over $\FP{q}$, as the following theorem shows.

\begin{thm}\label{T37}
Let $r,d\in \NN$ be such that $r\geq d$ and $\gcd(r,q-1)=1$. Let $H$ be a polynomial over $\FF{q}{2}$ of degree $d$ which permutes $\FP{q}$. Let $a,b\in \FF{q}{2}^\ast$ and $u,v\in \mu_{q+1}$ be such that $u\neq v$ and $\bigl(\frac{u}{v}\bigr)^d\neq \bigl(\frac{a}{b}\bigr)^{q-1}$. Let $w\in \FF{q}{2}^\ast$ be such that $w^{q-1}=\frac{u}{v}$. Then $F$ permutes $\FF{q}{2}$, where
\begin{equation}\label{Eq028}
F(x)=\begin{cases}
aw^{r-d}(v-u)^rx^r & \text{if }x^q-vx=0; \\
\bigl(x^q-vx\bigr)^r\left(aw^{-d}H\left(w\left(\frac{x^q-ux}{x^q-vx}\right)\right)+b\right) & \text{if }x^q-vx\neq 0.
\end{cases}
\end{equation}
In particular, if $\gcd(d,q-1)=1$ and $H(x)=x^d$, then $\widetilde{F}$ permutes $\FF{q}{2}$, where
\begin{equation}\label{Eq029}
\widetilde{F}(x)=\begin{cases}
aw^{r-d}(v-u)^rx^r & \text{if }x^q-vx=0; \\
a\bigl(x^q-vx\bigr)^{r-d}\bigl(x^q-ux\bigr)^d+b\bigl(x^q-vx\bigr)^r & \text{if }x^q-vx\neq 0.
\end{cases}
\end{equation}
\end{thm}

\noindent \textit{Proof}\quad Since \cref{Eq029} follows trivially from \cref{Eq028}, it suffices to verify \cref{Eq028}. Since $u,v\in \mu_{q+1}$, there exist $\omega_0,\omega_1\in \FF{q}{2}^\ast$ such that
\begin{align*}
u & =\omega_1^{q-1}; \\
v & =\omega_0^{q-1}.
\end{align*}
Since $u\neq v$, $W=\{\omega_0,\omega_1\}$ is $\FP{q}$-linearly independent, hence a basis for $\FF{q}{2}$ over $\FP{q}$. Since $w^{q-1}=\frac{u}{v}$, we may assume without loss of generality that
$$w=-\frac{\omega_1}{\omega_0}.$$
Then we define
\begin{equation}\label{Eq030}
\begin{aligned}
\gamma_0 & =aw^{-d}\omega_0^{-r}; \\
\gamma_1 & =b\omega_0^{-r}.
\end{aligned}
\end{equation}
By the initial assumptions,
$$\Bigl(\frac{u}{v}\Bigr)^d\neq \Bigl(\frac{a}{b}\Bigr)^{q-1}=w^{d(q-1)}\Bigl(\frac{\gamma_0}{\gamma_1}\Bigr)^{q-1}=\Bigl(\frac{u}{v}\Bigr)^d\Bigl(\frac{\gamma_0}{\gamma_1}\Bigr)^{q-1}.$$
Consequently, $\bigl(\frac{\gamma_0}{\gamma_1}\bigr)^{q-1}\neq 1$, meaning that $Y=\{\gamma_0,\gamma_1\}$ is $\FP{q}$-linearly independent, hence a basis for $\FF{q}{2}$ over $\FP{q}$. If $M_W$ is the Moore matrix associated with $W$, then the dual basis of $W$ is $B=\bigl\{\det({M_W})^{-1}\omega_0^q,-\det({M_W})^{-1}\omega_1^q\bigr\}$. Hence,
\begin{equation}\label{Eq031}
\begin{aligned}
Tr_q^{q^2}(\beta_0x) & =\det({M_W})^{-1}\left(\omega_1^qx-\omega_1x^q\right); \\
Tr_q^{q^2}(\beta_1x) & =\det({M_W})^{-1}\left(-\omega_0^qx+\omega_0x^q\right).
\end{aligned}
\end{equation}
We show that if $F$ is as defined by \cref{Eq028}, then $\det({M_W})^{-r}F$ permutes $\FF{q}{2}$. Clearly, $F(0)=0$. Now assume that $x\neq 0$. If $Tr_q^{q^2}(\beta_0x)\neq 0$ and $Tr_q^{q^2}(\beta_1x)=0$, then as per \cref{Eq031}, $x^q-vx=0$. Since $x\neq 0$, $x^{q-1}=v$. In \cref{Eq027}, if $H_0(x)=0$. Thus,
\begin{align*}
F_0(x)
& \stackrel{\mathmakebox[\widthof{=}]{(\ref{Eq027})}}{=} \Tr_q^{q^2}(\beta_0x)^r\gamma_0 \\
& \stackrel{\mathmakebox[\widthof{=}]{(\ref{Eq031})}}{=} \det({M_W})^{-r}\left(\omega_1^qx-\omega_1x^q\right)^r\gamma_0 \\
& =\det({M_W})^{-r}\gamma_0(-\omega_1)^r\bigl(x^{q-1}-\omega_1^{q-1}\bigr)^rx^r \\
& =\det({M_W})^{-r}\gamma_0\omega_0^r\Bigl(-\frac{\omega_1}{\omega_0}\Bigr)^r(v-u)^rx^r \\
& \stackrel{\mathmakebox[\widthof{=}]{(\ref{Eq030})}}{=} \det({M_W})^{-r}aw^{r-d}(v-u)^rx^r.
\end{align*}
If $Tr_q^{q^2}(\beta_1x)\neq 0$, then $x^q-vx\neq 0$. In \cref{Eq027}, $H_1(x)=H(x)$. Thus,
\begin{align*}
F_1(x)
& \stackrel{\mathmakebox[\widthof{=}]{(\ref{Eq027})}}{=} \Tr_q^{q^2}(\beta_1x)^r\biggl(H\biggl(\frac{\Tr_q^{q^2}(\beta_0x)}{\Tr_q^{q^2}(\beta_1x)}\biggr)\gamma_0+\gamma_1\biggr) \\
& \stackrel{\mathmakebox[\widthof{=}]{(\ref{Eq031})}}{=} \det({M_W})^{-r}\left(-\omega_0^qx+\omega_0x^q\right)^r\biggl(H\biggl(\frac{\omega_1^qx-\omega_1x^q}{-\omega_0^qx+\omega_0x^q}\biggr)\gamma_0+\gamma_1\biggr) \\
& \stackrel{\mathmakebox[\widthof{=}]{(\ref{Eq030})}}{=}\det({M_W})^{-r}\bigl(x^q-vx\bigr)^r\omega_0^r\biggl(aw^{-d}\omega_0^{-r}H\left(w\left(\frac{x^q-ux}{x^q-vx}\right)\right)+b\omega_0^{-r}\biggr) \\
& =\det({M_W})^{-r}\bigl(x^q-vx\bigr)^r\biggl(aw^{-d}H\left(w\left(\frac{x^q-ux}{x^q-vx}\right)\right)+b\biggr).
\end{align*}
So by \cref{T32}, $\det({M_W})^{-r}F$ is an $r$-homogeneous permutation polynomial of index $q+1$ over $\FF{q}{2}$. Therefore, the same is true for $F$.\mbox{}\hfill $\square$\\
\\
We give a small example of \cref{T36}, which has been verified by MAGMA.

\begin{exm}\label{E38}
Let $\omega$ be a primitive element in $\FF{8}{3}$. Then $W=\bigl\{1,\omega,\omega^2\bigr\}$ is a basis for $\FF{8}{3}$ over $\FP{8}$. If $B=\left\{\beta_0,\beta_1,\beta_2\right\}$ is the dual basis of $W$, then the polynomial $F$ permutes $\FF{q}{3}$, where
\[\resizebox{0.88\textwidth}{!}{
$\begin{aligned}
F(x)
& =
\left\{
\begin{array}{ll}
Tr_q^{q^2}(\beta_0x)^2 & \text{if }x\in \FP{8}; \\
Tr_q^{q^2}(\beta_0x)Tr_q^{q^2}(\beta_1x)+Tr_q^{q^2}(\beta_1x)^2\omega & \text{if }x\in \Spn_{\FP{8}}\left\{1,\omega\right\}\setminus \FP{8}; \\
Tr_q^{q^2}(\beta_0x)Tr_q^{q^2}(\beta_2x)+Tr_q^{q^2}(\beta_1x)Tr_q^{q^2}(\beta_2x)\omega+Tr_q^{q^2}(\beta_2x)^2\omega^2 & \text{if }x\in \FF{8}{3}\setminus \Spn_{\FP{8}}\left\{1,\omega\right\}.
\end{array}
\right.\\
& =\left\{
\begin{array}{ll}
\omega^{154}x^{128}+\omega^{147}x^{16}+\omega^{210}x^2 &\ \ \ \text{if }x\in \FP{8}; \\
\omega^{206}x^{128}+\omega^{367}x^{72}+\omega^{381}x^{65}+\omega^{283}x^{16}+\omega^{493}x^9+\omega^{464}x^2 &\ \ \ \text{if }x\in \Spn_{\FP{8}}\left\{1,\omega\right\}\setminus \FP{8}; \\
\omega^{433}x^{65} + \omega^{118}x^9+\omega^{398}x^2 &\ \ \ \text{if }x\in \FF{8}{3}\setminus \Spn_{\FP{8}}\left\{1,\omega\right\}.
\end{array}
\right.
\end{aligned}
$}\]
\end{exm}

For our next example of \cref{T32}, we need the following bijection over $\PG(n-1,q)$.

\begin{prp}\label{P39}
Let $q$ be an odd prime power, and let $\delta\in \FF{q}{n}^\ast$ be a primitive element in $\FP{q}^\ast$. Let $\alpha$ be any even power of $\delta$. Let $d\in \NN$ be coprime to $q-1$. Define $\widetilde{g}:\FP{q}^n\to \FP{q}^n$, where $\widetilde{g}\bigl(\vec{0}\bigr)=\vec{0}$, and for all $\lambda\in \FP{q}^\ast$,
\begin{enumerate}
    \item if $0\leq j\leq n-1$ and $j\neq 1$, then for all $x_0,\dots,x_{j-1}\in \FP{q}$,
    $$\lambda(x_0,\dots,x_{j-1},1,0,\dots,0)\stackrel{\widetilde{g}}{\mapsto}(x_0,\dots,x_{j-1},1,0,\dots,0);$$
    \item for all $x_0\in \FP{q}$,
    $$\lambda(x_0,1,\dots,0)\stackrel{\widetilde{g}}{\mapsto}
    \begin{cases}
    (x_0,1,0,\dots,0) & \text{if }x_0\text{ is an odd power of }\delta; \\
    \bigl(\alpha x_0^d,1,0,\dots,0\bigr) & \text{if }x_0\text{ is }0\text{ or an even power of }\delta.
    \end{cases}$$
\end{enumerate}
If $\widetilde{g}$ is written as $(g_0,\dots,g_{n-1})$, then $\overline{g}=(g_0:\dots:g_{n-1})$ permutes $\PG(n-1,q)$.
\end{prp}

\noindent \textit{Proof}\quad Since $\vec{0}$ is the only root of $\widetilde{g}$ in $\FP{q}^n$ and $\widetilde{g}\bigl(\lambda\vec{x}\bigr)=\widetilde{g}\bigl(\vec{x}\bigr)$ for $0\leq k\leq n-1$, $\overline{g}$ is well-defined over $\PG(n-1,q)$ as per \cref{L33}. Now we show that $\overline{g}$ permutes $\PG(n-1,q)$. For $0\leq j\leq n-1$, we recall that 
$$\overline{C_{q,n,j}}\stackrel{(\ref{Eq016})}{=}\{(x_0:\dots:x_{j-1}:1:0:\dots:0):x_0,\dots,x_{j-1}\in \FP{q}\}.$$
Since $\overline{g}$ fixes $\overline{C_{q,n,j}}$ for all $j\neq 1$ and $(0:1:0:\dots:0)$, it suffices to prove that $\overline{g}$ permutes $\overline{C_{q,n,1}}\setminus \{(0:1:0:\dots:0)\}$. Let $x_0,y_0$ be distinct elements in $\FP{q}$. If they are both odd powers of $\delta$, then
$$\overline{g}(x_0:1:0:\dots:0)=(x_0:1:0:\dots:0)\neq (y_0:1:0:\dots:0)=\overline{g}(y_0:1:0:\dots:0).$$
Similarly, if $x_0,y_0$ are both even powers of $\delta$, then since $\alpha\neq 0$ and $\gcd(d,q-1)=1$,
$$\overline{g}(x_0:1:0:\dots:0)=\bigl(\alpha x_0^d:1:0:\dots:0\bigr)\neq \bigl(\alpha y_0^d:1:0:\dots:0\bigr)=\overline{g}(y_0:1:0:\dots:0).$$
Hence, it remains to consider the case where $x_0$ is an odd power of $\delta$ and $y_0$ an even power of $\delta$. Assume towards a contradiction that $\overline{g}(x_0:1:0:\dots:0)=\overline{g}(y_0:1:0:\dots:0)$. Then $(x_0:1:0:\dots:0)=\bigl(\alpha y_0^d:1:0:\dots:0\bigr)$. Therefore, there exists a $\lambda\in \FP{q}^\ast$ such that $(x_0,1,0,\dots,0)=\lambda\bigl(\alpha y_0^d,1,0,\dots,0\bigr)$. Clearly, $\lambda=1$, and so $x_0=\alpha y_0^d$. Let $k,\ell,m\in \NN\sqcup \{0\}$ be such that $x_0=\delta^{2k+1},y_0=\delta^{2\ell}$ and $\alpha=\delta^{2m}$. Then $\delta^{2(k-d\ell-m)+1}=x_0\alpha^{-1}y_0^{-d}=1$, which is impossible since $2(k-d\ell-m)+1$ is odd but the multiplicative order of $\delta$ in $\FP{q}^\ast$ is $q-1$, which is even due to the assumption that $q$ is odd. Consequently, we must have that $\overline{g}(x_0:1:0:\dots:0)\neq \overline{g}(y_0:1:0:\dots:0)$. Therefore, $\overline{g}$ permutes $\PG(n-1,q)$.\mbox{}\hfill $\square$\\
\\
We also observe that if $W=\left\{\omega_j:0\leq j\leq n-1\right\}$ is any basis for $\FF{q}{n}$ over $\FP{q}$ with dual basis $B=\left\{\beta_j:0\leq j\leq n-1\right\}$, then for any  $x\in \FF{q}{n}^\ast$,
\begin{enumerate}
    \item $\Tr_q^{q^n}(\beta_1x)\neq 0$, $\Tr_q^{q^n}(\beta_ix)=0$ for all $i\geq 2$ and $\frac{\Tr_q^{q^n}(\beta_0x)}{\Tr_q^{q^n}(\beta_1x)}$ is an odd power of $\delta$ if and only if $x=\lambda\bigl(\delta^{2m+1}\omega_0+\omega_1\bigr)$ for some $\lambda\in \FP{q}^\ast$ and $m\in \NN\sqcup \{0\}$;
    \item $\Tr_q^{q^n}(\beta_1x)\neq 0$, $\Tr_q^{q^n}(\beta_ix)=0$ for all $i\geq 2$ and $\frac{\Tr_q^{q^n}(\beta_0x)}{\Tr_q^{q^n}(\beta_1x)}$ is an even power of $\delta$ if and only if $x=\lambda\bigl(\delta^{2m}\omega_0+\omega_1\bigr)$ for some $\lambda\in \FP{q}^\ast$ and $m\in \NN\sqcup \{0\}$.
\end{enumerate}
We can then construct homogeneous permutation polynomials of index $\frac{q^n-1}{q-1}$ over $\FF{q}{n}$ via \cref{T32,P39,Eq019} together with the above observations.

\begin{thm}\label{T310}
Let $q$ be an odd prime power, and let $\delta\in \FF{q}{n}^\ast$ be a primitive element in $\FP{q}^\ast$. Let $\alpha$ be any even power of $\delta$. Let $r,d\in \NN$ be coprime to $q-1$, and assume that $r\geq d$. Let $W=\left\{\omega_j:0\leq j\leq n-1\right\}$ and $Y=\left\{\gamma_j:0\leq j\leq n-1\right\}$ be any two bases for $\FF{q}{n}$ over $\FP{q}$. Let $B=\left\{\beta_j:0\leq j\leq n-1\right\}$ be the dual basis of $W$. Define $F$ piecewise, where
\begin{enumerate}
    \item if $j\neq 1$ and $x=\lambda(x_0\omega_0+\dots+x_{j-1}\omega_{j-1}+\omega_j)$ for some $\lambda\in \FP{q}^\ast$ and $x_0,\dots,x_{j-1}\in \FP{q}$, or if $j=1$ and $x=\lambda\bigl(\delta^{2m+1}\omega_0+\omega_1\bigr)$ for some $\lambda\in \FP{q}^\ast$ and $m\in \NN\sqcup \{0\}$, then
    $$F(x)=\sum_{k=0}^j\Tr_q^{q^n}(\beta_kx)\Tr_q^{q^n}(\beta_jx)^{r-1}\gamma_k;$$
    \item if $j=1$ and $x=\omega$ or $\lambda\bigl(\delta^{2m}\omega_0+\omega_1\bigr)$ for some $\lambda\in \FP{q}^\ast$ and $m\in \NN\sqcup \{0\}$, then
    $$F(x)=\alpha\Tr_q^{q^n}(\beta_0x)^d\Tr_q^{q^n}(\beta_1x)^{r-d}\gamma_0+\Tr_q^{q^n}(\beta_1x)^r\gamma_1.$$
\end{enumerate}
Then $F$ permutes $\FF{q}{n}$.
\end{thm}

We give a small example of \cref{T310}, which has been verified by MAGMA.

\begin{exm}\label{E311}
Let $\omega$ be a primitive element in $\FF{9}{3}$. Then $W=\bigl\{1,\omega,\omega^2\bigr\}$ is a basis for $\FF{9}{3}$ over $\FP{9}$, and $F$ permutes $\FF{9}{3}$, where
\[\resizebox{0.85\hsize}{!}{%
$F(x)=
\begin{cases}
\omega^{681}x^{81}+\omega^{561}x^9+\omega^{305}x & \text{if }x\in \FP{9}; \\
\omega^{309}x^{81}+\omega^{117}x^9+\omega^{553}x & \text{if }x=\lambda\bigl(\omega^{91m}+\omega\bigr), \lambda\in \FP{q}^\ast,m=1,3,5,7; \\
\omega^{615}x^{81}+\omega^{595}x^9+\omega^{670}x & \text{if }x=\omega\text{ or }\lambda\bigl(\omega^{91m}+\omega\bigr), \lambda\in \FP{q}^\ast,m=0,2,4,6; \\
\omega^{273}x & \text{otherwise}. \\
\end{cases}$%
}\]
\end{exm}

To conclude this section, we note that in all of our examples so far, the bijections $\overline{g}$ over $\PG(n-1,q)$ permute each of the sets $\overline{C_{q,n,j}}\ (0\leq j\leq n-1)$ given by \cref{Eq016}. However, for $\overline{g}$ to be a bijection over $\PG(n-1,q)$, it need not permute those sets. Hence, we present an example where $\overline{g}$ does not permute any of $\overline{C_{q,n,j}}\ (0\leq j\leq n-1)$.\\
\\
Let $\widetilde{g}:\FP{q}^3\to \FP{q}^3$, where $\widetilde{g}\bigl(\vec{0}\bigr)=\vec{0}$ and
\begin{align*}
\lambda(0,0,1)     & \stackrel{\widetilde{g}}{\mapsto}(1,0,0)     & & \text{for all }\lambda\in \FP{q}^\ast; \\
\lambda(1,0,0)     & \stackrel{\widetilde{g}}{\mapsto}(0,1,0)     & & \text{for all }\lambda\in \FP{q}^\ast; \\
\lambda(x_0,0,1)   & \stackrel{\widetilde{g}}{\mapsto}(x_0,1,0)   & & \text{for all }\lambda,x_0\in \FP{q}^\ast; \\
\lambda(x_0,1,0)   & \stackrel{\widetilde{g}}{\mapsto}(x_0,0,1)   & & \text{for all }\lambda\in \FP{q}^\ast\text{ and }x_0\in \FP{q}; \\
\lambda(x_0,x_1,1) & \stackrel{\widetilde{g}}{\mapsto}(x_0,x_1,1) & & \text{for all }\lambda,x_1\in \FP{q}^\ast\text{ and }x_0\in \FP{q}.
\end{align*}
Then it is easy to check that the corresponding $\overline{g}$ is a bijection over $\PG(2,q)$, where
\begin{align*}
\overline{g}(0:0:1)     & =(1:0:0) \\
\overline{g}(1:0:0)     & =(0:1:0); \\
\overline{g}(x_0:0:1)   & =(x_0:1:0)   & & \text{for all }x_0\in \FP{q}^\ast; \\
\overline{g}(x_0:1:0)   & =(x_0:0:1)   & & \text{for all }x_0\in \FP{q}; \\
\overline{g}(x_0:x_1:1) & =(x_0:x_1:1) & & \text{for all }x_0\in \FP{q}\text{ and }x_1\in \FP{q}^\ast.
\end{align*}

Using $\overline{g}$ together with \cref{T32}, we obtain the following result.

\begin{prp}\label{P312}
Let $W=\left\{\omega_j:0\leq j\leq 2\right\}$ and $Y=\left\{\gamma_j:0\leq j\leq 2\right\}$ be any two bases for $\FF{q}{3}$ over $\FP{q}$. Let $B=\left\{\beta_j:0\leq j\leq n-1\right\}$ be the dual basis of $W$. Let $r\in \NN$ be coprime to $q-1$. Then $F$ permutes $\FF{q}{3}$, where
\[
\resizebox{0.8\hsize}{!}{%
$F(x)=\left\{
\begin{array}{lll}
\Tr_q^{q^3}(\beta_2x)^r\gamma_0 & \text{if }x=x_2\omega_2\text{ and }x_2\in \FP{q}^\ast; \\
\\
\Tr_q^{q^3}(\beta_0x)^r\gamma_2 & \text{if }x=x_0\omega_0\text{ and }x_0\in \FP{q}; \\
\\
\Tr_q^{q^3}(\beta_0x)\Tr_q^{q^3}(\beta_2x)^{r-1}\gamma_0+\Tr_q^{q^3}(\beta_2x)^r\gamma_1 & \text{if }x=x_0\omega_0+x_2\omega_2\text{ and }x_0,x_2\in \FP{q}^\ast; \\
\\
\Tr_q^{q^3}(\beta_0x)\Tr_q^{q^3}(\beta_1x)^{r-1}\gamma_0+\Tr_q^{q^3}(\beta_1x)^r\gamma_2 & \text{if }x=x_0\omega_0+x_1\omega_1\text{ and }x_1\in \FP{q}^\ast; \\
\\
\Tr_q^{q^3}(\beta_0x)\Tr_q^{q^3}(\beta_2x)^{r-1}\gamma_0+ & \text{if }x=x_0\omega_0+x_1\omega_1+x_2\omega_2\text{ and }x_1,x_2\in \FP{q}^\ast. \\
\Tr_q^{q^3}(\beta_1x)\Tr_q^{q^3}(\beta_2x)^{r-1}\gamma_1+\Tr_q^{q^3}(\beta_2x)^r\gamma_2. 
\end{array}
\right.$%
}
\]
\end{prp}

We give a small example of \cref{P312}, which has been verified by MAGMA.

\begin{exm}\label{E313}
Let $\omega$ be a primitive element in $\FF{8}{3}$. Then $W=\bigl\{1,\omega,\omega^2\bigr\}$ is a basis for $\FF{8}{3}$ over $\FP{8}$, and $F$ permutes $\FF{8}{3}$, where
$$F(x)=\left\{
\begin{array}{lll}
\omega^{68}x^{64}+\omega^{264}x^8+\omega^{33}x  & \text{if }x=x_2\omega^2\text{ and }x_2\in \FP{q}^\ast; \\
\omega^{78}x^{64}+\omega^{330}x^8+\omega^{106}x & \text{if }x\in \FP{q}; \\
\omega^{87}x^{64}+\omega^{409}x^8+\omega^{303}x & \text{if }x=x_0+x_2\omega^2\text{ and }x_0,x_2\in \FP{q}^\ast;\\
\omega^{244}x^{64}+\omega^{64}x^8+\omega^{443}x & \text{if }x=x_0+x_1\omega\text{ and }x_0\in \FP{q},x_1\in \FP{q}^\ast; \\
\omega^{365}x & \text{if }x=x_0+x_1\omega+x_2\omega^2\text{ and }x_0\in \FP{q},x_1,x_2\in \FP{q}^\ast.   
\end{array}
\right.$$
\end{exm}

\section{A comparison with Hirschfeld's bijections}\label{sec4}
Having studied the GMTs, determined their inverses and presented some of their applications, a natural question to ask seems to be: are there other known classes of bijections defined on $\PG(n-1,q)$? The answer is affirmative. While studying the representation of the points and lines in projective spaces, Hirschfeld \cite{JH98a} introduced a class of bijections defined on $\PG(n-1,q)$. In this section, we compare the GMTs with Hirschfeld's mappings, and count the number of both classes of bijections. To describe Hirschfeld's results requires the notion of subprimitivity and projectivity.

\begin{dfn}[\cite{JH98a}]\label{D41}
A monic irreducible polynomial of degree $n$ over $\FP{q}$ is called subprimitive if the smallest $e\in \NN$ such that $f(x)\mid x^e-c$ for some $c\in \FP{q}^\ast$ is $\frac{q^n-1}{q-1}$. A root of a subprimitive polynomial is called a subprimitive root in $\FF{q}{n}$.
\end{dfn}

\begin{dfn}[\cite{JH98a}]\label{D42}
A projectivity is a bijection $\mathfrak{T}$ over $\PG(n-1,q)$ given by a non-singular matrix $M=M(\mathfrak{T})$ over $\FP{q}$, where $\mathfrak{T}(x_0:\dots:x_{n-1})=(x_0':\dots:x_{n-1}')$ if there exists a $\lambda\in \FP{q}^\ast$ such that $\lambda (x_0,\dots,x_{n-1})M=(x_0',\dots,x_{n-1}')$. A projectivity $\mathfrak{T}$ is cyclic if it permutes the elements of $\PG(n-1,q)$ in a single cycle. That is, if $(x_0:\dots:x_{n-1}),(y_0:\dots:y_{n-1})\in \PG(n-1,q)$, then there exists some $0\leq i\leq \frac{q^n-1}{q-1}-1$ such that $\mathfrak{T}^i(x_0:\dots:x_{n-1})=(y_0:\dots:y_{n-1})$ or, equivalently, $\lambda (x_0,\dots,x_{n-1})M^i=(y_0,\dots,y_{n-1})$ for some $\lambda\in \FP{q}^\ast$.
\end{dfn}

Let $\alpha$ be a subprimitive root in $\FF{q}{n}$. Let $f(x)=x^n-\displaystyle{\sum_{k=0}^{n-1} a_kx^k}\in \FP{q}[x]$ be the minimal polynomial of $\alpha$ over $\FP{q}$, and let
$$T=\begin{bmatrix}
0 & 1 & 0 & \dots & 0 \\
0 & 0 & 1 & \dots & 0 \\
\vdots & \vdots & \vdots & \ddots & \vdots \\
0 & 0 & 0 & \dots & 1 \\
a_0 & a_1 & a_2 & \dots & a_{n-1}
\end{bmatrix}$$
be the companion matrix of $f$. As per \cite[Theorem 4.2.]{JH98a}, a projectivity over $\PG(n-1,q)$ is cyclic if and only if its associated matrix is the companion matrix of a subprimitive polynomial of degree $n$ over $\FP{q}$. Thus, if $\mathfrak{T}$ is the projectivity over $\PG(n-1,q)$ given by \cref{D42} such that $M(\mathfrak{T})=T$, then $\mathfrak{T}$ is cyclic. So every element in $\PG(n-1,q)$ is the image of $(1:0:\dots:0)$ under $\mathfrak{T}^i$ for some $0\leq i\leq \frac{q^n-1}{q-1}-1$. For $0\leq i\leq \frac{q^n-1}{q-1}-1$, let $y_{\alpha,0}^{(i)},\dots,y_{\alpha,n-1}^{(i)}\in \FP{q}$ be the unique scalars such that
\begin{equation}\label{EqS41}
(1,0,\dots,0)T^i=\Bigl(y_{\alpha,0}^{(i)},\dots,y_{\alpha,n-1}^{(i)}\Bigr).
\end{equation}
Hirschfeld showed by induction that
\begin{equation}\label{EqS42}
\sum_{k=0}^{n-1}y_{\alpha,k}^{(i)}\alpha^k=\alpha^i\ \biggl(0\leq i \leq \frac{q^n-1}{q-1}-1\biggr).
\end{equation}
According to \cref{D41}, since $\alpha$ is a subprimitive root, the smallest $e\in \NN$ such that $\alpha^e\in \FP{q}^\ast$ is $\frac{q^n-1}{q-1}$. So \cref{EqS42} implies that if $0\leq i_1\neq i_2\leq \frac{q^n-1}{q-1}-1$, then there exists no $\lambda\in \FP{q}^\ast$ such that $\Bigl(y_{\alpha,0}^{(i_1)},\dots,y_{\alpha,n-1}^{(i_1)}\Bigr)=\lambda \Bigl(y_{\alpha,0}^{(i_2)},\dots,y_{\alpha,n-1}^{(i_2)}\Bigr)$. Hence, the equivalence classes $\Bigl(y_{\alpha,0}^{(i)}:\dots:y_{\alpha,n-1}^{(i)}\Bigr)$, where $0\leq i\leq \frac{q^n-1}{q-1}-1$, are pairwise distinct, meaning that
$$\PG(n-1,q)=\biggl\{\Bigl(y_{\alpha,0}^{(i)}:\dots:y_{\alpha,n-1}^{(i)}\Bigr):0\leq i\leq \frac{q^n-1}{q-1}-1\biggr\},$$
Let $P_\alpha=\bigl\{\alpha^i:0\leq i\leq \frac{q^n-1}{q-1}-1\bigr\}$. Hirschfeld concluded from the above observations that the mapping $\mathcal{H}_{\alpha}:\PG(n-1,q)\to P_\alpha$ is a bijection, where for $0\leq i\leq \frac{q^n-1}{q-1}-1$,
\begin{equation}\label{EqS43}
\mathcal{H}_{\alpha}\left(y_{\alpha,0}^{(i)}:\dots:y_{\alpha,n-1}^{(i)}\right)=\sum_{k=0}^{n-1}y_{\alpha,k}^{(i)}\alpha^k=\alpha^i.
\end{equation} 
The major differences between $\mathcal{H}_{\alpha}$ and the GMT $\psi_W$ given by \cref{Eq007} are as follows.
\begin{enumerate}
    \item First of all, $\psi_W(\PG(n-1,q))=\mu_{\frac{q^n-1}{q-1}}$, whereas $\mathcal{H}_{\alpha}(\PG(n-1,q))=P_\alpha$, where $\alpha$ is a subprimitive root in $\FF{q}{n}$. By \cref{D41}, the smallest $e\in \NN$ such that $\alpha^e=c$ for some $c\in \FP{q}^\ast$ is $e=\frac{q^n-1}{q-1}$. In particular, $P_\alpha=\mu_{\frac{q^n-1}{q-1}}$ if and only if $c=1$.
    \item Moreover, for $0\leq i\leq \frac{q^n-1}{q-1}-1$, $\mathcal{H}_{\alpha}^{-1}\left(\alpha^i\right)=\Bigl(y_{\alpha,0}^{(i)}:\dots:y_{\alpha,n-1}^{(i)}\Bigr)$, which is the equivalent class in $\PG(n-1,q)$ represented by the coordinate vector $\Bigl(y_{\alpha,0}^{(i)},\dots,y_{\alpha,n-1}^{(i)}\Bigr)$ of $\alpha^i$ with respect to the polynomial basis $\left\{\alpha^j:0\leq j\leq n-1\right\}$ for $\FF{q}{n}$ over $\FP{q}$. This is not the case for $\psi_W$, which is defined using an arbitrary basis $W=\left\{\omega_j:0\leq j\leq n-1\right\}$ for $\FF{q}{n}$ over $\FP{q}$. Indeed, if $x=(x_0\omega_0+\dots+x_{n-1}\omega_{n-1})^{q-1}\in \mu_{\frac{q^n-1}{q-1}}$, then $\psi_W^{-1}(x)=(x_0:\dots:x_{n-1})$, which is not necessarily the equivalence class in $\PG(n-1,q)$ represented by the coordinate vector of $x$ with respect to $W$ unless $q=2$, in which case $\mathcal{H}_{\alpha}=\psi_W$, where $W=\bigl\{\alpha^j:0\leq j\leq n-1\bigr\}$.
    \item The set of all GMTs from $\PG(n-1,q)$ to $\mu_{\frac{q^n-1}{q-1}}$ is closed under scalar multiplication by elements in $\mu_{\frac{q^n-1}{q-1}}$. Indeed, let $W=\left\{\omega_j:0\leq j\leq n-1\right\}$ be any basis for $\FF{q}{n}$ over $\FP{q}$, and let $c\in \mu_{\frac{q^n-1}{q-1}}$. Then $c=\widetilde{c}^{q-1}$ for some $\widetilde{c}\in \FF{q}{n}^\ast$. It is easy to show that $c\psi_W=\psi_{\widetilde{W}}$, where $\widetilde{W}=\widetilde{c}W=\left\{\widetilde{c}\omega_j:0\leq j\leq n-1\right\}$. However, if $\alpha,\beta,\gamma$ are subprimitive roots in $\FF{q}{n}$, then $\gamma^i\mathcal{H}_{\alpha}\neq \mathcal{H}_{\beta}$ for any $1\leq i\leq \frac{q^n-1}{q-1}-1$. Indeed, assume towards a contradiction that $\gamma^i\mathcal{H}_{\alpha}=\mathcal{H}_{\beta}$ for some $1\leq i\leq \frac{q^n-1}{q-1}-1$ and subprimitive roots $\alpha,\beta,\gamma$ in $\FF{q}{n}$. By \cref{EqS41}, $\Bigl(y_{\alpha,0}^{(0)},\dots,y_{\alpha,n-1}^{(0)}\Bigr)=\bigl(y_{\beta,0}^{(0)},\dots,y_{\beta,n-1}^{(0)}\bigr)=(1,0,\dots,0)$. Hence, $\gamma^i=\gamma^i\mathcal{H}_{\alpha}\Bigl(y_{\alpha,0}^{(0)},\dots,y_{\alpha,n-1}^{(0)}\Bigr)=\mathcal{H}_{\beta}\Bigl(y_{\beta,0}^{(0)},\dots,y_{\beta,n-1}^{(0)}\Bigr)=\beta^0=1$. Since $\gamma$ is a subprimitive root, the smallest $e\in \NN$ such that $\gamma^e=c$ for some $c\in \FP{q}^\ast$ is $\frac{q^n-1}{q-1}$. Hence, $\frac{q^n-1}{q-1}\mid i$, which is impossible since $1\leq i\leq \frac{q^n-1}{q-1}-1$.
\end{enumerate}

Moreover, there are more GMTs than Hirschfeld's bijections over $\PG(n-1,q)$, as the following result shows.

\begin{prp}\label{P43}
Over $\PG(n-1,q)$, let $M(n,q)$ and $H(n,q)$ denote the number of GMTs and that of Hirschfeld's bijections, respectively. Then
$$\frac{M(n,q)}{H(n,q)}>q^{n-1}\prod_{i=1}^{n-2}\left(q^n-q^i\right).$$
\end{prp}

\noindent \textit{Proof}\quad Let $W=\left\{\omega_j:0\leq j\leq n-1\right\}$ be a basis for $\FF{q}{n}$ over $\FP{q}$. Let $\sigma\in S_n$, and define $W^{\sigma}=\left\{\theta_j:0\leq j\leq n-1\right\}$, where $\theta_j=\omega_{\sigma(j)}$ for each $j$. Then $\psi_W=\psi_{W^\sigma}$ if and only if $\sigma=\text{id}$. Indeed, if $\sigma=\text{id}$, then clearly $\psi_W=\psi_{W^\sigma}$. Conversely, assume that $\psi_W=\psi_{W^\sigma}$. Let $0\leq j_0\leq n-1$. If $(x_0:\dots:x_{n-1})\in \PG(n-1,q)$ is such that $x_{j_0}=1$ and $x_j=0$ for all $j\neq j_0$, then $\psi_W(x_0:\dots:x_{n-1})=\omega_{j_0}^{q-1}$, and $\psi_{W^{\sigma}}(x_0:\dots:x_{n-1})=\theta_{j_0}^{q-1}=\omega_{\sigma(j_0)}^{q-1}$. Since $\psi_W=\psi_{W^\sigma}$, $\omega_{j_0}^{q-1}=\omega_{\sigma(j_0)}^{q-1}$. Since $W$ is $\FP{q}$-linearly independent, $\omega_{j_0}=\omega_{\sigma(j_0)}$, and so $j_0=\sigma(j_0)$. Since this is true for all $j_0$, $\sigma=\text{id}$. Meanwhile, for any two bases $W=\left\{\omega_j:0\leq j\leq n-1\right\},Y=\left\{\gamma_j:0\leq j\leq n-1\right\}$ for $\FF{q}{n}$ over $\FP{q}$, $\psi_W=\psi_Y$ if and only if $W=\lambda Y=\left\{\lambda \omega:\omega\in Y\right\}$ for some $\lambda\in \FP{q}^\ast$. Indeed, if the latter condition is true, then $\psi_W=\psi_Y$. Conversely, assume that $\psi_W=\psi_Y$. Since $\psi_W(1:0:\dots:0)=\psi_Y(1:0:\dots:0),\dots,\psi_W(0:0:\dots:1)=\psi_Y(0:0:\dots:1)$, we know that for each $0\leq j\leq n-1$, there exists some $\lambda_j\in \FP{q}^\ast$ such that $\omega_j=\lambda_j\gamma_j$. Since $\psi_W(1:1:\dots:1)=\psi_Y(1:1:\dots:1)$, there exists some $\lambda\in \FP{q}^\ast$ such that $$\lambda_0\gamma_0+\dots+\lambda_{n-1}\gamma_{n-1}=\omega_0+\dots+\omega_{n-1}=\lambda(\gamma_0+\dots+\gamma_{n-1}).$$
Since $Y$ is $\FP{q}$-linearly independent, $\lambda_0=\dots=\lambda_{n-1}=\lambda$. Thus, $M(n,q)$ is the number of \textit{ordered} bases for $\FF{q}{n}$ over $\FP{q}$ divided by $q-1$. That is,
$$M(n,q)=\frac{\left(q^n-1\right)\left(q^n-q\right)\dots \left(q^n-q^{n-1}\right)}{q-1}.$$
Finally, we determine $H(n,q)$. If $\alpha_1,\alpha_2$ are distinct subprimitive roots $\FF{q}{n}$, then $\mathcal{H}_{\alpha_1}\neq \mathcal{H}_{\alpha_2}$. Indeed, by \cref{EqS41}, $(1,0,\dots,0)T=(0,1,0,\dots,0)$ for all companion matrices $T$. So for $j=1,2$, $\mathcal{H}_{\alpha_j}(0:1:\dots:0)=\alpha_j$ due to \cref{EqS43}. In particular, $\mathcal{H}_{\alpha_1}\neq \mathcal{H}_{\alpha_2}$. Consequently, $H(n,q)$ is the number of subprimitive roots in $\FF{q}{n}$, i.e., $n$ times the number of subprimitive polynomials of degree $n$ over $\FP{q}$. By \cite[Eq. (1.12)]{JH98a}, there are $\bigl(\frac{q-1}{n}\bigr)\varphi\bigl(\frac{q^n-1}{q-1}\bigr)$ subprimitive polynomials of degree $n$ over $\FP{q}$, where $\varphi$ is the Euler totient function. Thus,
$$H(n,q)=(q-1)\varphi\bigl(\frac{q^n-1}{q-1}\bigr),$$
In particular, we know that
$$\frac{M(n,q)}{H(n,q)}>\frac{M(n,q)}{q^n-1}=q^{n-1}\prod_{i=1}^{n-2}\left(q^n-q^i\right),$$
which completes the proof.\mbox{}\hfill $\square$

\section{Conclusions}\label{sec5}
In this paper, we generalized the M\"obius transformation to a class of bijections from $\PG(n-1,q)$ to $\mu_{\frac{q^n-1}{q-1}}$ for an arbitrary positive integer $n\geq 2$, and determined the inverses of the GMTs. As their applications, we constructed various classes of homogeneous permutation polynomials of index $\frac{q^n-1}{q-1}$ over $\FF{q}{n}$. We also showed that the GMTs are different from Hirschfeld's bijections  from $\PG(n-1,q)$ to the $i$-th powers of subprimitive roots, where $0\leq i\leq \frac{q^n-1}{q-1}-1$, thereby extending the family of known bijections defined on $\PG(n-1,q)$ with good properties and applications.

\setcounter{secnumdepth}{0}

\end{document}